\documentclass[12pt]{article}

\usepackage{amsthm,amsmath,amssymb,mathtools}
\usepackage{ascmac}
\usepackage{color}
\usepackage{enumitem}
\usepackage{hyperref}
\usepackage{tikz}
\usepackage{subcaption}

\mathtoolsset{showonlyrefs=true}

\usepackage{sty/macros}
\usepackage{sty/mathoperators}
\usepackage{sty/theorems}

\hypersetup{
     colorlinks=true,
     linkcolor=black,
     filecolor=blue,
     citecolor = black,
     urlcolor=cyan,
     }

\def\@themcountersep{}
\usepackage{fancybox,color}
\definecolor{lred}{rgb}{1,0.8,0.5}
\definecolor{lblue}{rgb}{0.8,0.8,1}
\definecolor{dred}{rgb}{0.6,0,0}
\definecolor{dblue}{rgb}{0,0,0.7}
\definecolor{violet}{rgb}{0.5804,0.0000,0.8275}
\definecolor{purple}{rgb}{0.2400,0.5700,0.2500}
\definecolor{TGreen}{rgb}{0,0.50,0.10}
\title{Exact SDP relaxations for quadratic programs with bipartite graph structures}

\makeatletter
\let\@fnsymbol\@arabic
\makeatother

\author{
\normalsize
    Godai Azuma\thanks{Department of  Mathematical and Computing Science,
        Tokyo Institute of Technology, 2-12-1-W8-29 Oh-Okayama, Meguro-ku, Tokyo 152-8552, Japan.
        ({\tt azuma.g.aa@m.titech.ac.jp, mituhiro@is.titech.ac.jp},  {\tt Makoto.Yamashita@c.titech.ac.jp}).
        The research of Makoto Yamashita was partially supported by JSPS KAKENHI Grant Number JP20H04145.}
\and
\normalsize
        Mituhiro Fukuda\footnotemark[1]\;\textsuperscript{,}\thinspace
        \thanks{Department of Computer Science, Institute of Mathematics and Statistics, University of S\~ao Paulo, Rua do Mat\~ao, 1010, Cidade Universit\'aria, S\~ao Paulo, SP, 05508-090, Brazil, and currently at S\~ao Paulo State Technological College, Praia Grande, Pra\chige{c}a 19 de Janeiro, 144, Praia Grande, SP, 11700-100, Brazil. The research of Mituhiro Fukuda was supported by grants 2020/04585-7 and 2018/24293-0 from the S\~ao Paulo Research Foundation (FAPESP).}
\and
\normalsize
	Sunyoung Kim\thanks{Department of Mathematics, Ewha W. University, 52 Ewhayeodae-gil, Sudaemoon-gu,
	Seoul 	03760, Korea  ({\tt skim@ewha.ac.kr}). This work was supported
        by  NRF 2021-R1A2C1003810.}
\and
\normalsize
        Makoto Yamashita\footnotemark[1]
        }

\begin{document}
\maketitle

\begin{abstract} \noindent
For nonconvex quadratically constrained quadratic programs (QCQPs),
we first show that, under certain feasibility conditions, the standard semidefinite (SDP) relaxation is exact for QCQPs with bipartite graph structures.
The exact optimal solutions are obtained by examining the dual SDP relaxation
and the rank of the optimal solution of this dual SDP relaxation under  strong duality.
Our  results on the QCQPs generalize the results on
QCQP with sign-definite bipartite graph structures, QCQPs with forest structures,
and QCQPs with nonpositive off-diagonal data elements.
Second, we propose a conversion method from QCQPs with no particular structure
to the ones with bipartite graph structures. 
As a result, we demonstrate that a wider class of QCQPs can be exactly solved by the SDP relaxation.
Numerical instances are presented for illustration.
\end{abstract}

\vspace{0.5cm}

\noindent
{\bf Key words. } Quadratically constrained quadratic programs, Exact semidefinite relaxations,
Bipartite graph, Sign-indefinite QCQPs, Rank of aggregated sparsity matrix.

\vspace{0.5cm}

\noindent
{\bf AMS Classification. }
90C20,  	
90C22,  	
90C25, 	
90C26.  	

\section{Introduction} \label{sec:introduction}

We consider nonconvex quadratically constrained quadratic programs (QCQPs) of the form
\begin{align}
    \label{eq:hqcqp} \tag{$\PC$}
	\begin{array}{rl}
        \min & \trans{\x}Q^0 \x \\
        \subto & \trans{\x}Q^p \x \leq b_p, \quad p \in [m],
    \end{array}
\end{align}
where $Q^0, \ldots, Q^m \in \SymMat^n$, $\b \in \Real^m$, $\x \in \Real^n$, and $[m]$ denotes
the set $\left\{i \in \Natural\,\middle|\, 1 \leq i \leq m \right\}$.
We use $\SymMat^n$ to denote
the space of $n \times n$ symmetric matrices.
A general form of QCQPs with linear terms
\begin{align*}
\begin{array}{rl}
\min & \trans{\x}Q^0 \x + \trans{(\q^0)}\x   \\
\subto & \trans{\x}Q^p \x + \trans{(\q^p)}\x \leq b_p  \quad p \in [m],
    \end{array}
\end{align*}
 can be represented in the form of \eqref{eq:hqcqp}
using a new variable $x_0$ such that $x_0^2 = 1$, where $\q^0, \ldots, \q^m \in \Real^n$.
For simplicity, we describe QCQPs as  \eqref{eq:hqcqp}
and we assume that \eqref{eq:hqcqp} is feasible in this paper.

Nonconvex QCQPs \eqref{eq:hqcqp} are known to be NP-hard in general,  however, finding the exact solution of
some class of QCQPs has been  a popular subject \cite{Azuma2021,Burer2019,Jeyakumar2014,kim2003exact,Sojoudi2014exactness,Wang2021geometric,Wang2021tightness}
 as they can provide  solutions for important applications
 formulated as QCQPs \eqref{eq:hqcqp}. They include optimal power flow problems~\cite{Lavaei2012,Zhou2019}, pooling problems~\cite{kimizuka2019solving},
 sensor network localization problems~\cite{BISWAS2004,KIM2009,SO2007},
quadratic assignment problems~\cite{PRendl09,ZHAO1998},
the max-cut problem~\cite{Geomans1995}.
Moreover, it is well-known that polynomial optimization problems can be recast as QCQPs.

By replacing $\x\trans{\x}$ with a rank-1 matrix $X \in \SymMat^n$ in \eqref{eq:hqcqp}
and removing the rank constraint of $X$,
the standard (Shor) SDP relaxation
and its dual problem can be expressed as
\begin{align}
    \label{eq:hsdr} \tag{$\PC_R$} &
    \begin{array}{rl}
        \min & \ip{Q^0}{X} \\
        \subto & \ip{Q^p}{X} \le b_p, \quad p \in [m], \\
               & X \succeq O,
    \end{array} \\
    \label{eq:hsdrd} \tag{$\DC_R$} &
    \begin{array}{rl}
        \max & \trans{-\b}\y \\
        \subto & S(\y) := Q^0 + \sum\limits_{p=1}^m y_p Q^p  \succeq O,
        \quad \y \ge \0,
    \end{array}
\end{align}
where $\ip{Q^p}{X}$ denotes the Frobenius inner product of $Q^p$ and $X$, i.e., $\ip{Q^p}{X} \coloneqq  \sum_{i,j} Q^p_{ij} X_{ij}$,
and $X \succeq O$ means that $X$ is positive semidefinite.
The SDP relaxation provides a lower bound of the optimal value of \eqref{eq:hqcqp} in general.
When the SDP relaxation ({$\PC_R$}) provides a rank-1  solution $X$, we say that
the SDP relaxation is exact. In this case, the exact optimal solution
and exact optimal value can be computed in polynomial time.
A second-order cone programming (SOCP) relaxation can be obtained by
further relaxing the positive semidefinite constraint $X \succeq O$, for instance,
requiring  all $2 \times 2$ principal submatrices of $X$ to be positive
semidefinite~\cite{kim2003exact, sheen2020exploiting}.
For QCQPs with a certain sparsity structure, e.g., forest structures,
the SDP relaxation coincides with the SOCP relaxation.

In this paper, we present a wider class of QCQPs that can be solved exactly with the SDP relaxation by extending
the results in  \cite{Azuma2021} and \cite{Sojoudi2014exactness}.  The extension is based on that trees or forests are  bipartite graphs
and that QCQPs with no structure and the same sign of $Q_{ij}^p$ for $ p=0,1, \ldots,m$ can be transformed into
ones with bipartite structures.
Sufficient conditions  for the exact
 SDP relaxation of QCQP \eqref{eq:hqcqp} are described. These conditions are called exactness conditions in the subsequent discussion.
We mention that our results on the exact SDP relaxation  is obtained  by investigating the rank of $S(\y)$ in the dual of SDP relaxation ({$\DC_R$}).

When discussing the exact optimal solution  of nonconvex QCQPs, convex relaxations of QCQPs such as the SDP or
SOCP have played a pivotal role.
In particular,
the signs of the elements in the data matrices $Q^0, \ldots, Q^m$
as in \cite{kim2003exact,Sojoudi2014exactness} and
graph structures such as forests \cite{Azuma2021} and bipartite structures \cite{Sojoudi2014exactness} have been used to identify the classes
of nonconvex QCQPs whose exact optimal solution can be attained via the SDP relaxation.
QCQPs with nonpositive off-diagonal data matrices were shown to have an exact SDP and SOCP relaxation~\cite{kim2003exact}.
This result was generalized by Sojoudi and Lavaei~\cite{Sojoudi2014exactness}
with  a sufficient condition  that can be tested by
the sign-definiteness based on the cycles in the aggregated sparsity pattern graph induced from the nonzero elements of data matrices in \eqref{eq:hqcqp}.
A finite set $\{Q^0_{ij}, Q^1_{ij}, \ldots, Q^m_{ij}\} \subseteq \Real$ is called sign-definite
if the elements of the set are either all nonnnegative or all nonpositive.
We note that these results are obtained by analyzing the primal problem ($\PC_R$).
 For general QCQPs with no particular structure,
Burer and Ye in   \cite{Burer2019} presented  sufficient conditions for the exact semidefinite formulation
with a polynomial-time checkable polyhedral system.
From the dual SDP relaxation \eqref{eq:hsdrd} using strong duality,
they proposed
an LP-based technique to detect the exactness of the SDP relaxation of QCQPs
consisting of diagonal matrices $Q^0, \ldots, Q^m$ and linear terms.
Azuma et al.~\cite{Azuma2021} presented  related results on QCQPs with forest structures.

With respect to the exactness conditions,
Yakubovich's S-lemma~\cite{Polik2007, Yakubovich1971}
(also known as S-procedure) can be regarded as one of the most important results.
It showed that the trust-region subproblem,
a subclass of QCQPs with only one constraint ($m = 1$) and $Q^1 \succeq O$, always admits an exact SDP relaxation.
Under some mild assumptions,
Wang and Xia~\cite{Wang2015} generalized this result to QCQPs with two constraints ($m = 2$)
and any matrices satisfing $Q^1 = -Q^2$ but not necessarily being positive semidefinite.
For the extended trust-region subproblem
whose constraints consist of one ellipsoid and linear inequalities,
 the exact SDP relaxation has been studied by
Jeyakumar and Li~\cite{Jeyakumar2014}. They proved that
the SDP relaxation of the extended trust-region subproblem is exact if
the algebraic multiplicity of the minimum eigenvalue of $Q^0$ is strictly greater than
the dimension of the space spanned by the coefficient vectors of the linear inequalities.
This condition was slightly improved by Hsia and Sheu~\cite{Hsia2013}.
In addition, Locatelli~\cite{Locatelli2016} introduced
a new exactness condition for the extended trust-region subproblem
based on the KKT conditions 
and proved that it is more general than the previous results.

%

A different approach on  the exactness of the SDP relaxation for QCQPs is to study  the convex hull exactness, i.e., the coincidence of
 the convex hull of the epigraph of a QCQP  and the projected epigraph of its SDP relaxation.
Wang and K{\dotlessi}l{\dotlessi}n\chige{c}-Karzan in~\cite{Wang2021tightness}
presented sufficient conditions for the convex hull exactness under the condition that
the feasible set $\Gamma \coloneqq \{\y \geq \0 \,|\, S(\y) \succeq O\}$ of \eqref{eq:hsdrd} is polyhedral.
Their results were improved in~\cite{Wang2021geometric} by eliminating this condition.
The rank-one generated (ROG) property, a geometric property,
was employed by Argue et al.~\cite{Argue2020necessary}
to evaluate the feasible set of the SDP relaxation.
In their paper, they proposed sufficient conditions that
the feasible set of the SDP relaxation is ROG, and
connected the ROG property with both the objective value and the exactness of the convex hull.

We describe our contributions:
\begin{itemize} \vspace{-2mm}
\item
We first show that
if the aggregated sparsity pattern graph is connected and bipartite
and a feasibility checking system constructed from  QCQP \eqref{eq:hqcqp} is infeasible,
then the  SDP relaxation is exact  in section~\ref{sec:main}.
It is a  polynomial-time  method as the systems can be represented as SDPs.
This result can be regarded as an extension of Azuma et al.~\cite{Azuma2021}
 in the sense that
 the aggregated sparsity pattern 
was generalized from  forests to  bipartite. We should mention that the signs of  data are irrelavant.
We give in section~\ref{sec:example} two numerical examples of QCQPs which can be shown to 
have exact SDP relaxations by our method, but fails to meet the conditions for real-valued
QCQP of \cite{Sojoudi2014exactness}.

\item
We propose a conversion method to  derive a   bipartite graph structure in  \eqref{eq:hqcqp} from QCQPs with no apparent  structure,
so that the SDP relaxation of the resulting QCQP provides the exact optimal solution.
More precisely, for every off-diagonal index $(i,j)$, if the set $\{Q^0_{ij},\ldots,Q^m_{ij}\}$ is sign-definite,
i.e., either all nonnegative or all nonpositive,
then any QCQP \eqref{eq:hqcqp} can be transformed into  nonnegative off-diagonal QCQPs
with bipartite aggregated sparsity
by introducing a new variable $\z \coloneqq -\x$ and a new constraint $\|\x + \z\|_2^2 \leq 0$, which
  covers a result for the real-valued QCQP proposed in \cite{Sojoudi2014exactness}.

\item We also show that the known results on the exactness of QCQPs
where (a) all the off-diagonal elements are sign-definite and the aggregated sparsity pattern graph is forest
or (b) all the off-diagonal elements are nonpositive can be proved using our method.

\item  For disconnected pattern graphs,
a perturbation of the objective function  is introduced, as  in~\cite{Azuma2021}, in section~\ref{sec:perturbation}
to demonstrate that a QCQP is exact
if there exists a sequence of perturbed problems converging to the QCQP
while maintaining the exactness of their SDP relaxation
 under assumptions weaker than \cite{Azuma2021}.

\end{itemize}

Throughout this paper, the following example is used to illustrate the difference between our result and previous works.
\begin{example}
    \label{examp:cycle-graph-4-vertices}
    \begin{align*}
       \begin{array}{rl}
       \min & \trans{\x} Q^0 \x \\
       \subto & \trans{\x} Q^1 \x \leq 10, \quad \trans{\x} Q^2 \x \leq 10,
         \end{array}
  \end{align*}
  where \begin{align*}
        & Q^0 = \begin{bmatrix}
             0 & -2 & 0 &  2 \\ -2 &  0 & -1 &  0 \\
             0 & -1 & 5 &  1 \\  2 &  0 &  1 & -4 \end{bmatrix},  \ 
          Q^1 = \begin{bmatrix}
             5 &  2 & 0 &  1 \\  2 & -1 &  3 &  0 \\
             0 &  3 & 3 & -1 \\  1 &  0 & -1 &  4 \end{bmatrix}, \  
         Q^2 = \begin{bmatrix}
            -1 &  1 & 0 &  0 \\  1 &  4 & -1 &  0 \\
             0 & -1 & 6 &  1 \\  0 &  0 &  1 & -2 \end{bmatrix}.
    \end{align*}
\end{example}

\noindent
Although Example \ref{examp:cycle-graph-4-vertices} does not satisfy
 the sign-definiteness,
the proposed method can successfully show that  the SDP relaxation is exact.

The rest of this paper is organized as follows.
In section~\ref{sec:preliminaries},
 the aggregated sparsity pattern of QCQPs and the sign-definiteness are defined
and related works on the exactness of the SDP relaxation for QCQPs
with some aggregated sparsity pattern are described.
Sections~\ref{sec:main} and \ref{sec:perturbation} include the main results of this paper.
In section~\ref{sec:main},  the assumptions necessary for the exact SDP relaxation are described,
and sufficient conditions for the exact  SDP relaxation are presented
under  the connectivity of the aggregated sparsity pattern.
In section~\ref{sec:perturbation},
we show that the sufficient conditions can be extended to
QCQPs which do not satisfy the connectivity condition.
The perturbation results on the exactness are utilized to remove the connectivity condition.
In section~\ref{sec:example}, we also provide specific numerical instances to  compare our result with the existing work and
 illustrate our method.
Finally, we conclude in section~\ref{sec:conclution}.

\section{Preliminaries} \label{sec:preliminaries}

We denote
the $n$-dimensional Euclidean space by $\Real^n$
and the nonnegative orthant of $\Real^n$ by $\Real_+^n$.
We  write the zero vector and the vector of all ones as $\0 \in \Real^n$ and  $\1 \in \Real^n$,  respectively.
We also write $M \succeq O$ and $M \succ O$ to indicate that
the matrix $M$ is positive semidefinite and positive definite, respectively.
We use $[n]:=\left\{i \in \Natural\,\middle|\, 1 \leq i \leq n\right\}$
and $[n, m]:=\left\{i \in \Natural\,\middle|\, n \leq i \leq m\right\}$.
The graph $G(\VC, \EC)$ denotes an undirected graph with
the vertex set $\VC$ and the edge set $\EC$.
We sometimes write $G$ if the vertex and edge sets are clear.


\subsection{Aggregated sparsity pattern}
The aggregated sparsity pattern of the SDP relaxation, defined from the data matrices $Q^p \, (p \in [0, m])$, is used
to describe the sparsity structure of QCQPs.
Let $\VC = [n]$ denote the set of indices of
rows and columns of $n \times n$ symmetric matrices.
Then, the set of indices
\begin{equation*}
    \EC = \left\{
        (i, j) \in \VC \times \VC \,\middle|\,
        \text{$i \neq j$ and $Q^p_{ij} \ne 0$ for some $p \in [0, m]$}
    \right\}
\end{equation*}
is called the aggregated sparsity pattern
for both a given QCQP \eqref{eq:hqcqp} and its SDP relaxation  \eqref{eq:hsdr}.
If $\EC$ denotes the set of edges of a graph with vertices $\VC$,
the graph $G(\VC, \EC)$ is called the aggregated sparsity pattern graph.
If $\EC$ corresponds to an adjacent matrix $\QC$ of $n$ vertices,
$\QC$ is called the aggregated sparsity pattern matrix.

\begin{figure}[t]
    \centering
    \begin{minipage}{0.30\textwidth}
        \tikzset{every picture/.style={line width=0.75pt}}
        \begin{tikzpicture}[x=0.75pt,y=0.75pt,yscale=-0.6,xscale=0.6]
            \draw   (140,40) -- (40,140) -- (140,140) -- (40,40) -- cycle ;
            \draw  [fill={rgb, 255:red, 255; green, 255; blue, 255 }  ,fill opacity=1 ] (10,40) .. controls (10,23.43) and (23.43,10) .. (40,10) .. controls (56.57,10) and (70,23.43) .. (70,40) .. controls (70,56.57) and (56.57,70) .. (40,70) .. controls (23.43,70) and (10,56.57) .. (10,40) -- cycle ;
            \draw  [fill={rgb, 255:red, 255; green, 255; blue, 255 }  ,fill opacity=1 ] (10,140) .. controls (10,123.43) and (23.43,110) .. (40,110) .. controls (56.57,110) and (70,123.43) .. (70,140) .. controls (70,156.57) and (56.57,170) .. (40,170) .. controls (23.43,170) and (10,156.57) .. (10,140) -- cycle ;
            \draw  [fill={rgb, 255:red, 255; green, 255; blue, 255 }  ,fill opacity=1 ] (110,140) .. controls (110,123.43) and (123.43,110) .. (140,110) .. controls (156.57,110) and (170,123.43) .. (170,140) .. controls (170,156.57) and (156.57,170) .. (140,170) .. controls (123.43,170) and (110,156.57) .. (110,140) -- cycle ;
            \draw  [fill={rgb, 255:red, 255; green, 255; blue, 255 }  ,fill opacity=1 ] (110,40) .. controls (110,23.43) and (123.43,10) .. (140,10) .. controls (156.57,10) and (170,23.43) .. (170,40) .. controls (170,56.57) and (156.57,70) .. (140,70) .. controls (123.43,70) and (110,56.57) .. (110,40) -- cycle ;

            \draw (40,40) node  [font=\large] [align=left] {\begin{minipage}[lt]{40.8pt}\setlength\topsep{0pt}
            \begin{center}
            {\fontfamily{pcr}\selectfont 1}
            \end{center}\end{minipage}};
            \draw (140,40) node  [font=\large] [align=left] {\begin{minipage}[lt]{40.8pt}\setlength\topsep{0pt}
            \begin{center}
            {\fontfamily{pcr}\selectfont 2}
            \end{center}\end{minipage}};
            \draw (140,140) node  [font=\large] [align=left] {\begin{minipage}[lt]{40.8pt}\setlength\topsep{0pt}
            \begin{center}
            {\fontfamily{pcr}\selectfont 4}
            \end{center}\end{minipage}};
            \draw (40,140) node  [font=\large] [align=left] {\begin{minipage}[lt]{40.8pt}\setlength\topsep{0pt}
            \begin{center}
            {\fontfamily{pcr}\selectfont 3}
            \end{center}\end{minipage}};
        \end{tikzpicture}
    \end{minipage}
    \begin{minipage}{0.30\textwidth}
        \begin{equation*}
            \QC = \begin{bmatrix}
                \star & \star & 0 & \star \\ \star & \star & \star & 0 \\
                0 & \star & \star & \star \\ \star & 0 & \star & \star
                \end{bmatrix}.
        \end{equation*}
    \end{minipage}
    \caption{The aggregated sparsity pattern graph and matrix of Example~\ref{examp:cycle-graph-4-vertices}.  $\star$ denotes an arbitrary  value.}
    \label{fig:example-aggregated-sparsity}
\end{figure}
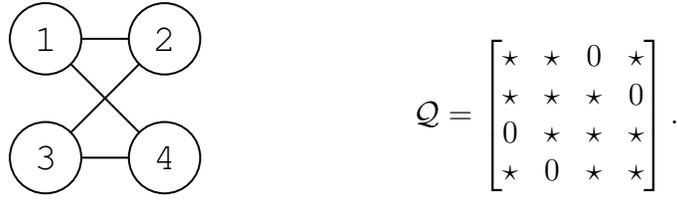

Consider the QCQP in Example \ref{examp:cycle-graph-4-vertices} as an illustrative example.
%
%
\noindent
As $(1, 3)$th and $(2, 4)$th elements are zeros in $Q^0, Q^1, Q^2$,
the aggregated sparsity pattern graph is a cycle  with 4 vertices as
shown in Figure~\ref{fig:example-aggregated-sparsity}. We see that
the graph  has only one cycle with 4 vertices.
This graph is the simplest of connected bipartite graphs with cycles.

For the discussion on QCQPs with sign-definiteness, we adopt the following notation
from \cite{Sojoudi2014exactness}.
We define the sign  $\sigma_{ij}$  of  each edge in $\VC \times \VC$ as
\begin{equation} \label{eq:definition-sigma-ij}
    \sigma_{ij} = \begin{cases}
        \quad +1 \quad & \text{if $Q^0_{ij}, \ldots, Q^m_{ij} \geq 0$,} \\
        \quad -1 \quad & \text{if $Q^0_{ij}, \ldots, Q^m_{ij} \leq 0$,} \\
        \quad  0 \quad & \text{otherwise.}
    \end{cases}
\end{equation}
Obviously,
$\sigma_{ij} \in \{-1, +1\}$ if and only if $\{Q^0_{ij}, \ldots, Q^m_{ij}\}$ is sign-definite.

Sojoudi and Lavaei~\cite{Sojoudi2014exactness}
proposed the following condition for exactness.
\vspace{-2mm}
\begin{theorem}[{\cite[Theorem 2]{Sojoudi2014exactness}}] \label{thm:sojoudi-theorem}
    The SOCP relaxation and the SDP relaxation of \eqref{eq:hqcqp} are exact
    if both of the following hold:
    \begin{align}
        &\sigma_{ij} \neq 0, && \forall (i, j) \in \EC, \label{eq:sign-constraint-sign-definite} \\
        \prod_{(i,j) \in \mathcal{C}_r} & \sigma_{ij} = (-1)^{\left|\mathcal{C}_r\right|}, && \forall r \in \{1,\ldots, \kappa\}, \label{eq:sign-constraint-simple-cycle}
    \end{align}
    where the set of cycles $\mathcal{C}_1, \ldots, \mathcal{C}_\kappa \subseteq \EC$ denotes a cycle basis for $G$.
\end{theorem}
\vspace{-2mm}
With the aggregated sparsity pattern graph $G$ of a given QCQP,
they presented the following corollary:
\begin{coro}[{\cite[Corollary 1]{Sojoudi2014exactness}}]
    \label{coro:sojoudi-corollary1}
    The SDP relaxation and the SOCP relaxation of \eqref{eq:hqcqp} are exact
    if one of the following holds:
    \begin{enumerate}[label=(\alph*)]
        \item $G$ is forest with $\sigma_{ij} \in \{-1, 1\}$ for all $(i, j) \in \EC$, \label{cond:sojoudi-forest}
        \item $G$ is bipartite with $\sigma_{ij} = 1$ for all $(i, j) \in \EC$, \label{cond:sojoudi-bipartite}
        \item $G$ is arbitrary with $\sigma_{ij} = -1$ for all $(i, j) \in \EC$. \label{cond:sojoudi-arbitrary}
    \end{enumerate}
\end{coro}

\subsection{Conditions for exact SDP relaxations with forest structures}
Recently,
Azuma et al.~\cite{Azuma2021} proposed
a method to decide the exactness of the SDP relaxation of QCQPs with forest structures.
The forest-structured QCQPs or their SDP relaxation have no cycles in their aggregated sparsity pattern graph.
In their work,
the rank of the dual SDP relaxation was determined using feasibility systems under the following assumption:
\begin{assum}
    \label{assum:previous-assumption}
    The following conditions hold for \eqref{eq:hqcqp}:
    \begin{enumerate}[label=(\roman*)]
        \item there exists $\bar{\y} \geq 0$ such that $\sum \bar{y}_p Q^p \succ O$, and \label{assum:previous-assumption-1}
        \item \eqref{eq:hsdr} has an interior feasible point. \label{assum:previous-assumption-2}
    \end{enumerate}
\end{assum}
\noindent
We note that Assumption  \ref{assum:previous-assumption} is used to derive 
 strong duality of the SDP relaxation
and the boundedness of the feasible set.
More precisely, for $\bar{\y}$ in Assumption~\ref{assum:previous-assumption},
multiplying $\ip{Q^p}{X} \leq b_p$ by $\bar{y}_p$ and adding together leads to
\begin{equation*}
    \ip{\left(\sum_{p = 1}^m \bar{y}_pQ^p\right)}{X} \leq \trans{\b}\bar{\y},
\end{equation*}
which implies that the feasible set of $X$ is bounded from $X \succeq O$.

We describe the result in \cite{Azuma2021} for our subsequent discussion.
\begin{prop}[\cite{Azuma2021}] \label{prop:forest-results}
    Assume that a given QCQP satisfies Assumption~\ref{assum:previous-assumption},
    and that its aggregated sparsity pattern graph $G(\VC, \EC)$ is a forest.
    The problem \eqref{eq:hsdr} is exact
    if, for all $(k, \ell) \in \EC$, the following system has no solutions:
    \begin{equation} \label{eq:system-zero}
        \y \geq 0,\; S(\y) \succeq O,\; S(\y)_{k\ell} = 0.
    \end{equation}
\end{prop}
\noindent
The above feasibility system,  formulated as SDPs, can be checked in polynomial time
since the number of edges of a forest graph with $n$ vertices is at most $n - 1$.

\section{Conditions for exact SDP relaxations with connected bipartite  structures} \label{sec:main}
Throughout this section,
we assume that the aggregated sparsity pattern graph $G(\VC, \EC)$ of a QCQP is connected and bipartite.
Under this assumption,
we present sufficient conditions for the SDP relaxation to be exact.
The main result described in Theorem~\ref{thm:system-based-condition-connected} in this section is extended to
the ones for the disconnected aggregated sparsity in section~\ref{sec:perturbation}.

Assumption~\ref{assum:previous-assumption}
has been introduced
only to derive the strong duality which is used in the proof of Proposition~\ref{prop:forest-results}.
Instead of Assumption~\ref{assum:previous-assumption}, we introduce
Assumption~\ref{assum:new-assumption}.
In  Remark~\ref{rema:comparison-assumption} below, we will consider a relation between
Assumptions~\ref{assum:previous-assumption} and
\ref{assum:new-assumption}.
\begin{assum} \label{assum:new-assumption}
    The following two conditions hold:
    \begin{enumerate}[label=(\roman*)]
        \item \label{assum:new-assumption-1}
        the sets of optimal solutions for \eqref{eq:hsdr} and \eqref{eq:hsdrd} are nonempty; and
        \item \label{assum:new-assumption-2}
        at least one of the following two conditions holds:
        \begin{enumerate}[label=(\alph*)]
            \item \label{assum:new-assumption-2-1}
                the feasible set of \eqref{eq:hsdr} is bounded; or
            \item \label{assum:new-assumption-2-2}
                the set of optimal solutions for \eqref{eq:hsdrd} is bounded.
        \end{enumerate}
    \end{enumerate}
\end{assum}
\noindent
The following lemma states that  strong duality holds under Assumption~\ref{assum:new-assumption}. 

\begin{lemma}
 \label{lem:feasible-set-strong-duality}
  If Assumption~\ref{assum:new-assumption} is satisfied,
  strong duality holds between \eqref{eq:hsdr} and \eqref{eq:hsdrd}, that is,
\eqref{eq:hsdr}   and \eqref{eq:hsdrd} have optimal solutions
and their optimal values are finite and equal.
\end{lemma}
\begin{proof}
    Since either the set of optimal solutions for \eqref{eq:hsdr} or
    that for \eqref{eq:hsdrd} is nonempty and bounded,
    Corollary~4.4 of Kim and Kojima~\cite{kim2021strong} indicates that
    the optimal values of \eqref{eq:hsdr} and \eqref{eq:hsdrd} are finite and equal.
\end{proof}

\begin{rema} \label{rema:comparison-assumption}
    Assumption~\ref{assum:new-assumption} is weaker than Assumption~\ref{assum:previous-assumption}.
    To compare these assumptions, we suppose that there exists $\bar{\y} \geq 0$ such that $\sum_p \bar{y}_pQ^p \succ O$.
    Then, there obviously exists sufficiently large $\lambda > 0$
    such that
    \begin{equation*}
        \lambda \bar{\y} \geq \0 \quad \text{and} \quad Q^0 + \sum_p \lambda\bar{y}_pQ^p \succ O,
    \end{equation*}
    which implies \eqref{eq:hsdrd} has an interior feasible point.
    It follows that the set of optimal solutions of \eqref{eq:hsdr} is bounded.
    Similarly, since \eqref{eq:hsdr} has an interior point by Assumption~\ref{assum:previous-assumption},
	the set of optimal solutions of \eqref{eq:hsdrd} is also bounded.
	This indicates Assumption~\ref{assum:new-assumption} {\it \ref{assum:new-assumption-1}} and
	{\it \ref{assum:new-assumption-2}\ref{assum:new-assumption-2-2}}.

    In addition,
    as mentioned right after Assumption~\ref{assum:previous-assumption},
    the feasible set of \eqref{eq:hsdr} is bounded.
    Thus,
    Assumption~\ref{assum:new-assumption}
    {\it \ref{assum:new-assumption-2}\ref{assum:new-assumption-2-1}} is also satisfied,
	under Assumption~\ref{assum:previous-assumption}.
\end{rema}

\subsection{Bipartite sparsity pattern matrix} \label{ssec:bipartite-matrix}
For a given matrix $M \in \SymMat^n$,
a sparsity pattern graph $G(\VC, \EC_M)$ can be defined by the vertex set and edge set:
\begin{equation*}
    \VC = [n], \quad
    \EC_M = \left\{(i, j) \in \VC \times \VC \,\middle|\, M_{ij} \neq 0\right\}.
\end{equation*}
Conversely, if $(i, j) \not\in \EC_M$, then the $(i, j)$th element of $M$ must be zero.

The graph $G(\VC, \EC)$ is called bipartite if
its vertices can be divided into two disjoint sets $\LC$ and $\RC$ such that
no two vertices in the same set are adjacent.
Equivalently, a bipartite $G$ is a graph with no odd cycles.
If $G(\VC, \EC)$ is bipartite, it can be represented with $G(\LC, \RC, \EC)$,
where $\LC$ and $\RC$ are disjoint sets of vertices.
The sets $\LC$ and $\RC$ are sometimes called parts of the bipartite graph $G$.

The following lemma is an immediate consequence of Proposition 1 of \cite{grone1992nonchordal}.
It shows that the rank of a nonnegative positive semidefinite matrix can be bounded below by $n - 1$
under some sparsity conditions if the sum of every row of the matrix is positive. 
We utilize Lemma~\ref{lemma:bipartite-rank} to estimate the rank of solutions of the dual SDP relaxation,
and establish conditions for the exact SDP relaxation in this section. 
\begin{lemma}[{\cite[Proposition 1]{grone1992nonchordal}}] \label{lemma:bipartite-rank}
    Let $M \in \Real^{n \times n}$ be a nonnegative and positive semidefinite matrix with $M\1 > \0$.
    If the sparsity pattern graph of $M$ is bipartite and connected,
    then $\rank(M) \geq n - 1$.
\end{lemma}

As the aggregated sparsity pattern graph $G$ composed from
$Q_0, Q_1, \dots, Q_m$ is used
to investigate the exactness of the SDP relaxation of a QCQP,
 the sparsity pattern graph of the matrix $S(\y)$ in the dual of the SDP relaxation is clearly a subgraph of $G$.
As a result, if $G$ is bipartite,
then the rank of $S(\y)$ can be estimated by Lemma~\ref{lemma:bipartite-rank}
since $S(\y)$ is also bipartite.
This will be used in the proof of Theorem~\ref{thm:system-based-condition-connected}.

\subsection{Main results} \label{ssec:main-connected}
We present our main results, that is, 
sufficient conditions for the SDP relaxation of the QCQP with bipartite structures to be exact.
%
\begin{theorem}
    \label{thm:system-based-condition-connected}
    Suppose that Assumption~\ref{assum:new-assumption} holds
    and the aggregated sparsity pattern $G(\VC, \EC)$ is a bipartite graph.
    Then, \eqref{eq:hsdr} is exact if
    \begin{itemize}
        \item $G(\VC, \EC)$ is connected,
        \item
            for all $(k, \ell) \in \EC$, the following system has no solutions:
            \begin{equation} \label{eq:system-nonpositive}
                \y \geq \0,\; S(\y) \succeq O,\; S(\y)_{k\ell} \leq 0.
            \end{equation}
    \end{itemize}
\end{theorem}
\begin{proof}
    Let $X^*$ be any optimal solution for \eqref{eq:hsdr} which exists by Assumption~\ref{assum:new-assumption}.
    By Lemma~\ref{lem:feasible-set-strong-duality},
    the optimal values of \eqref{eq:hsdr} and \eqref{eq:hsdrd} are finite and equal.
    Thus, there exists an optimal solution $\y^*$ for \eqref{eq:hsdrd}
    such that the complementary slackness holds, i.e.,
    \begin{equation*}
        X^* S(\y^*) = O.
    \end{equation*}
    Since $\y^* \geq \0$ and $S(\y^*) \succeq O$,
    by the infeasibility of \eqref{eq:system-nonpositive},
    we obtain $S(\y^*)_{k\ell} > 0$ for every $(k, \ell) \in \EC$.
    Furthermore, for each $i \in \VC$, the $i$th element of $S(\y^*)\1$ is
    \begin{equation*}
        [S(\y^*)\1]_i
        = \sum_{j = 1}^n S(\y^*)_{ij}
        = S(\y^*)_{ii} + \sum_{(i,j) \in \EC} S(\y^*)_{ij}
        > 0.
    \end{equation*}
    By Lemma~\ref{lemma:bipartite-rank},
    $\rank\left\{S(\y^*)\right\} \geq n - 1$.
    From the Sylvester’s rank inequality~\cite{Anton2014},
    \begin{equation*}
        \rank(X^*)
        \leq n - \rank\left\{S(\y^*)\right\} + \rank\left\{X^*S(\y^*)\right\}
        \leq n - (n - 1)
        = 1.
    \end{equation*}
    Therefore, the SDP relaxation is exact.
\end{proof}


The exactness of a given QCQP can be determined
by checking the infeasibility of $|\EC|$ systems.
Since \eqref{eq:system-nonpositive} can be formulated as
an SDP with the objective function $0$,
 checking their infeasibility is not difficult.

Compared with Proposition~\ref{prop:forest-results} in \cite{Azuma2021},
Theorem~\ref{thm:system-based-condition-connected} can determine the exactness of a wider class of QCQPs
in terms of the required assumption and sparsity.
As mentioned in Remark~\ref{rema:comparison-assumption},
the  assumptions in Theorem \ref{thm:system-based-condition-connected}  are weaker
than those in  Proposition~\ref{prop:forest-results},
and the aggregated sparsity pattern of $G$ is extended from forest graphs to bipartite graphs.


\subsection{Nonnegative off-diagonal QCQPs} \label{ssec:nonnegative-offdiagonal-connected}

We can also prove a known result by Theorem~\ref{thm:system-based-condition-connected},
i.e.,
the exactness of the SDP relaxation for QCQPs with nonnegative off-diagonal data matrices $Q^0, \ldots, Q^m$,
which was
referred as Corollary~\ref{coro:sojoudi-corollary1}\ref{cond:sojoudi-bipartite} above and was proved in~\cite{Sojoudi2014exactness}.
The aggregated sparsity pattern graph $G(\VC, \EC)$ is assumed to be connected
and $Q^0_{ij} > 0$ for all $(i, j) \in \EC$ in this subsection.
These assumptions will be relaxed in
section~\ref{ssec:nonnegative-offdiagonal}.

\begin{coro} \label{coro:nonnegative-offdiagonal-connected}
    Suppose that Assumption~\ref{assum:new-assumption} holds,
    and the aggregated sparsity pattern graph $G(\VC, \EC)$ of \eqref{eq:hqcqp}
    is bipartite and connected.
    If $Q^0_{ij} > 0$ for all $(i, j) \in \EC$ and
    $Q^p_{ij} \geq 0$ for all $(i, j) \in \EC$ and all $p \in [m]$,
    then the SDP relaxation is exact.
\end{coro}
\begin{proof}
    Let $\hat{\y} \ge \0$ be any nonnegative vector satisfying $S(\hat{\y}) \succeq O$.
    By the assumption, for any $(i, j) \in \EC$,
    \begin{equation*}
        S(\hat{\y})_{ij}
        = Q^0_{ij} + \sum_{p \in [m]} \hat{y}_p Q^p_{ij}
        \geq Q^0_{ij}
        > 0.
    \end{equation*}
    Hence, the system \eqref{eq:system-nonpositive} for every $(i, j) \in \EC$ has no solutions.
    Therefore, by Theorem~\ref{thm:system-based-condition-connected},
    the SDP relaxation is exact.
\end{proof}

\subsection{Conversion to QCQPs with bipartite structures} \label{ssec:comparison}

We show that a QCQP can be transformed into an equivalent QCQP with bipartite structures.
We then compare Theorem~\ref{thm:system-based-condition-connected}
with Theorem~\ref{thm:sojoudi-theorem}. 
As our result has been obtained by
the rank of the dual SDP \eqref{eq:hsdrd} via  strong duality while
the result in \cite{Sojoudi2014exactness} is from the evaluation of \eqref{eq:hsdr},  the classes of QCQPs that can be solved exactly with the SDP relaxation become different.
In this section, we show that
a class of QCQPs obtained by Theorem~\ref{thm:system-based-condition-connected}
under Assumption~\ref{assum:new-assumption} is wider than those by Theorem~\ref{thm:sojoudi-theorem}.

To transform a QCQP into an equivalent QCQP with bipartite structures and to apply Theorem~\ref{thm:system-based-condition-connected},
we define a diagonal matrix $D^p \in \SymMat^n$ with a positive number
from the diagonal of $Q^p$ for every $p$.
In addition, off-diagonal elements of $Q^p$ are divided into
two nonnegative symmetric matrices $2N^p_+,\, 2N^p_- \in \SymMat^n$ according to their signs
such that $Q^p = D^p + 2N^p_+ - 2N^p_-$.
More precisely, for an arbitrary positive number $\delta > 0$,
\begin{align*}
    D^p_{ii} &= Q^p_{ii} + 2 \delta, \\
    2[N^p_+]_{ij} &= \begin{cases}
       + Q^p_{ij} & \text{if $i \neq j$ and $Q^p_{ij} > 0$}, \\
        0 & \text{otherwise,}
    \end{cases} \\
    2[N^p_-]_{ij} &= \begin{cases}
      - Q^p_{ij} & \text{if $i \neq j$ and $Q^p_{ij} < 0$}, \\
      2 \delta & \text{if $i = j$}, \\
        0 & \text{otherwise.}
    \end{cases}
\end{align*}
We introduce a new variable $\z$ such that $\z \coloneqq -\x$. Then,
\begin{equation*}
    \trans{\x}Q^p\x =
        \trans{\begin{bmatrix} \x \\ \z \end{bmatrix}}
        \begin{bmatrix} D^p + 2N^p_+ & N^p_- \\ N^p_- & O \end{bmatrix}
        \begin{bmatrix} \x \\ \z \end{bmatrix},
\end{equation*}
The constraint $\z = -\x$ can be expressed as
$\|\x + \z\|^2 \leq 0$,
which can be written as
\begin{equation*}
    \trans{(\x + \z)}(\x + \z)
    = \trans{\begin{bmatrix} \x \\ \z \end{bmatrix}}
        \begin{bmatrix} I & I \\ I & I \end{bmatrix}
        \begin{bmatrix} \x \\ \z \end{bmatrix}
    \leq 0.
\end{equation*}
Thus, we have an equivalent QCQP:
\begin{equation} \label{eq:decomposed-hqcqp}
    \begin{array}{rl}
        \min & \trans{\begin{bmatrix} \x \\ \z \end{bmatrix}}
            \begin{bmatrix} D^0 + 2N^0_+ & N^0_- \\ N^0_- & O \end{bmatrix}
            \begin{bmatrix} \x \\ \z \end{bmatrix} \\
        \subto & \trans{\begin{bmatrix} \x \\ \z \end{bmatrix}}
            \begin{bmatrix} D^p + 2N^p_+ & N^p_- \\ N^p_- & O \end{bmatrix}
            \begin{bmatrix} \x \\ \z \end{bmatrix} \leq b_p, \quad p \in [m], \\
        & \trans{\begin{bmatrix} \x \\ \z \end{bmatrix}}
            \begin{bmatrix} I & I \\ I & I \end{bmatrix}
            \begin{bmatrix} \x \\ \z \end{bmatrix} \leq 0.
    \end{array}
\end{equation}
Note that
\eqref{eq:decomposed-hqcqp} includes $m + 1$ constraints and
all off-diagonal elements of data matrices are nonnegative since $N^p_+$ and $N^p_-$ are nonnegative.
Let $\bar{G}(\bar{\VC}, \bar{\EC})$ denote
the aggregated sparsity pattern graph of \eqref{eq:decomposed-hqcqp}.
The number of vertices in $\bar{G}$ is twice as many as that 
in $G$ due to the additional variable $\z$.
If $\bar{G}$ is bipartite and $Q^0_{ij} \neq 0$ for all $(i, j) \in \EC$,
the SDP relaxation of \eqref{eq:decomposed-hqcqp} is exact
since the assumptions of Corollary~\ref{coro:nonnegative-offdiagonal-connected} are satisfied.

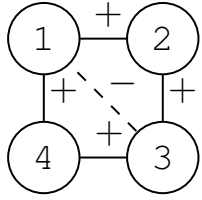
\begin{figure}[t]
    \centering
    \tikzset{every picture/.style={line width=0.75pt}} 
    \begin{tikzpicture}[x=0.75pt,y=0.75pt,yscale=-0.6,xscale=0.6]
        \draw [dash pattern={on 4.5pt off 4.5pt}, line width=0.75]    (40,40) -- (140,140) ;
        \draw   (40,40) -- (140,40) -- (140,140) -- (40,140) -- cycle ;
        \draw  [fill={rgb, 255:red, 255; green, 255; blue, 255 }  ,fill opacity=1 ] (10,40) .. controls (10,23.43) and (23.43,10) .. (40,10) .. controls (56.57,10) and (70,23.43) .. (70,40) .. controls (70,56.57) and (56.57,70) .. (40,70) .. controls (23.43,70) and (10,56.57) .. (10,40) -- cycle ;
        \draw  [fill={rgb, 255:red, 255; green, 255; blue, 255 }  ,fill opacity=1 ] (10,140) .. controls (10,123.43) and (23.43,110) .. (40,110) .. controls (56.57,110) and (70,123.43) .. (70,140) .. controls (70,156.57) and (56.57,170) .. (40,170) .. controls (23.43,170) and (10,156.57) .. (10,140) -- cycle ;
        \draw  [fill={rgb, 255:red, 255; green, 255; blue, 255 }  ,fill opacity=1 ] (110,140) .. controls (110,123.43) and (123.43,110) .. (140,110) .. controls (156.57,110) and (170,123.43) .. (170,140) .. controls (170,156.57) and (156.57,170) .. (140,170) .. controls (123.43,170) and (110,156.57) .. (110,140) -- cycle ;
        \draw  [fill={rgb, 255:red, 255; green, 255; blue, 255 }  ,fill opacity=1 ] (110,40) .. controls (110,23.43) and (123.43,10) .. (140,10) .. controls (156.57,10) and (170,23.43) .. (170,40) .. controls (170,56.57) and (156.57,70) .. (140,70) .. controls (123.43,70) and (110,56.57) .. (110,40) -- cycle ;

        \draw (40,40) node  [font=\large] [align=left] {\begin{minipage}[lt]{40.8pt}\setlength\topsep{0pt}
        \begin{center}
        {\fontfamily{pcr}\selectfont 1}
        \end{center}\end{minipage}};

        \draw (140,40) node  [font=\large] [align=left] {\begin{minipage}[lt]{40.8pt}\setlength\topsep{0pt}
        \begin{center}
        {\fontfamily{pcr}\selectfont 2}
        \end{center}\end{minipage}};

        \draw (140,140) node  [font=\large] [align=left] {\begin{minipage}[lt]{40.8pt}\setlength\topsep{0pt}
        \begin{center}
        {\fontfamily{pcr}\selectfont 3}
        \end{center}\end{minipage}};

        \draw (40,140) node  [font=\large] [align=left] {\begin{minipage}[lt]{40.8pt}\setlength\topsep{0pt}
        \begin{center}
        {\fontfamily{pcr}\selectfont 4}
        \end{center}\end{minipage}};

        \draw (91.33,64.33) node [anchor=north west][inner sep=0.75pt]  [font=\large] [align=left] {$\displaystyle -$};
        \draw (79.67,7) node [anchor=north west][inner sep=0.75pt]  [font=\large] [align=left] {$\displaystyle +$};
        \draw (80.33,106.67) node [anchor=north west][inner sep=0.75pt]  [font=\large] [align=left] {$\displaystyle +$};
        \draw (42,71) node [anchor=north west][inner sep=0.75pt]  [font=\large] [align=left] {$\displaystyle +$};
        \draw (142,71) node [anchor=north west][inner sep=0.75pt]  [font=\large] [align=left] {$\displaystyle +$};
    \end{tikzpicture}
    \caption{
        An aggregated sparsity pattern graph with edge signs.
        The solid and dashed lines show that
        the corresponding $\sigma_{ij}$ are $+1$ and $-1$, respectively.
        Both lines indicate the existence of nonzero elements in some $Q^p$.
        }
    \label{fig:example-edge-signs}
\end{figure}

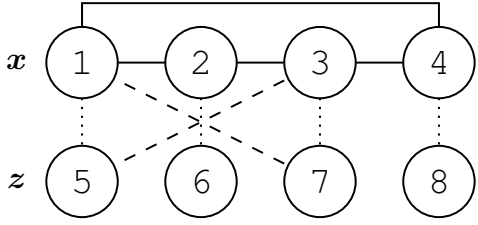
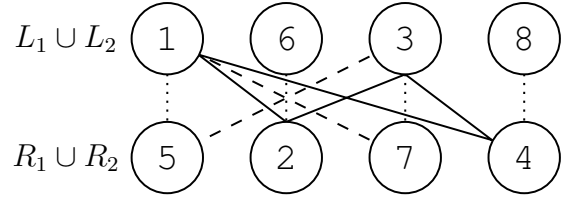
\begin{figure}[t]
    \centering
    \begin{subfigure}[t]{0.45\textwidth}
        \centering
        \tikzset{every picture/.style={line width=0.75pt}} 
        \begin{tikzpicture}[x=0.75pt,y=0.75pt,yscale=-0.6,xscale=0.6]
            \draw  [dash pattern={on 4.5pt off 4.5pt}]  (70,60) -- (270,160) ;
            \draw  [dash pattern={on 4.5pt off 4.5pt}]  (70,160) -- (270,60) ;
            \draw  [dash pattern={on 0.84pt off 2.51pt}]  (270,60) -- (270,160) ;
            \draw  [dash pattern={on 0.84pt off 2.51pt}]  (70,60) -- (70,160) ;
            \draw  [dash pattern={on 0.84pt off 2.51pt}]  (170,60) -- (170,160) ;
            \draw  [dash pattern={on 0.84pt off 2.51pt}]  (370,60) -- (370,160) ;
            \draw   (370,10) -- (70,10) -- (70,60) -- (370,60) -- cycle ;

            \draw  [fill={rgb, 255:red, 255; green, 255; blue, 255 }  ,fill opacity=1 ] (40,60) .. controls (40,43.43) and (53.43,30) .. (70,30) .. controls (86.57,30) and (100,43.43) .. (100,60) .. controls (100,76.57) and (86.57,90) .. (70,90) .. controls (53.43,90) and (40,76.57) .. (40,60) -- cycle ;
            \draw  [fill={rgb, 255:red, 255; green, 255; blue, 255 }  ,fill opacity=1 ] (240,60) .. controls (240,43.43) and (253.43,30) .. (270,30) .. controls (286.57,30) and (300,43.43) .. (300,60) .. controls (300,76.57) and (286.57,90) .. (270,90) .. controls (253.43,90) and (240,76.57) .. (240,60) -- cycle ;
            \draw  [fill={rgb, 255:red, 255; green, 255; blue, 255 }  ,fill opacity=1 ] (140,60) .. controls (140,43.43) and (153.43,30) .. (170,30) .. controls (186.57,30) and (200,43.43) .. (200,60) .. controls (200,76.57) and (186.57,90) .. (170,90) .. controls (153.43,90) and (140,76.57) .. (140,60) -- cycle ;
            \draw  [fill={rgb, 255:red, 255; green, 255; blue, 255 }  ,fill opacity=1 ] (340,60) .. controls (340,43.43) and (353.43,30) .. (370,30) .. controls (386.57,30) and (400,43.43) .. (400,60) .. controls (400,76.57) and (386.57,90) .. (370,90) .. controls (353.43,90) and (340,76.57) .. (340,60) -- cycle ;
            \draw  [fill={rgb, 255:red, 255; green, 255; blue, 255 }  ,fill opacity=1 ] (40,160) .. controls (40,143.43) and (53.43,130) .. (70,130) .. controls (86.57,130) and (100,143.43) .. (100,160) .. controls (100,176.57) and (86.57,190) .. (70,190) .. controls (53.43,190) and (40,176.57) .. (40,160) -- cycle ;
            \draw  [fill={rgb, 255:red, 255; green, 255; blue, 255 }  ,fill opacity=1 ] (240,160) .. controls (240,143.43) and (253.43,130) .. (270,130) .. controls (286.57,130) and (300,143.43) .. (300,160) .. controls (300,176.57) and (286.57,190) .. (270,190) .. controls (253.43,190) and (240,176.57) .. (240,160) -- cycle ;
            \draw  [fill={rgb, 255:red, 255; green, 255; blue, 255 }  ,fill opacity=1 ] (140,160) .. controls (140,143.43) and (153.43,130) .. (170,130) .. controls (186.57,130) and (200,143.43) .. (200,160) .. controls (200,176.57) and (186.57,190) .. (170,190) .. controls (153.43,190) and (140,176.57) .. (140,160) -- cycle ;
            \draw  [fill={rgb, 255:red, 255; green, 255; blue, 255 }  ,fill opacity=1 ] (340,160) .. controls (340,143.43) and (353.43,130) .. (370,130) .. controls (386.57,130) and (400,143.43) .. (400,160) .. controls (400,176.57) and (386.57,190) .. (370,190) .. controls (353.43,190) and (340,176.57) .. (340,160) -- cycle ;

            \draw (170,60) node  [font=\large] [align=left] {\begin{minipage}[lt]{40.8pt}\setlength\topsep{0pt}
            \begin{center}
            {\fontfamily{pcr}\selectfont 2}
            \end{center}\end{minipage}};

            \draw (370,60) node  [font=\large] [align=left] {\begin{minipage}[lt]{40.8pt}\setlength\topsep{0pt}
            \begin{center}
            {\fontfamily{pcr}\selectfont 4}
            \end{center}\end{minipage}};

            \draw (70,60) node  [font=\large] [align=left] {\begin{minipage}[lt]{40.8pt}\setlength\topsep{0pt}
            \begin{center}
            {\fontfamily{pcr}\selectfont 1}
            \end{center}\end{minipage}};

            \draw (270,60) node  [font=\large] [align=left] {\begin{minipage}[lt]{40.8pt}\setlength\topsep{0pt}
            \begin{center}
            {\fontfamily{pcr}\selectfont 3}
            \end{center}\end{minipage}};

            \draw (170,160) node  [font=\large] [align=left] {\begin{minipage}[lt]{40.8pt}\setlength\topsep{0pt}
            \begin{center}
            {\fontfamily{pcr}\selectfont 6}
            \end{center}\end{minipage}};

            \draw (370,160) node  [font=\large] [align=left] {\begin{minipage}[lt]{40.8pt}\setlength\topsep{0pt}
            \begin{center}
            {\fontfamily{pcr}\selectfont 8}
            \end{center}\end{minipage}};

            \draw (70,160) node  [font=\large] [align=left] {\begin{minipage}[lt]{40.8pt}\setlength\topsep{0pt}
            \begin{center}
            {\fontfamily{pcr}\selectfont 5}
            \end{center}\end{minipage}};

            \draw (270,160) node  [font=\large] [align=left] {\begin{minipage}[lt]{40.8pt}\setlength\topsep{0pt}
            \begin{center}
            {\fontfamily{pcr}\selectfont 7}
            \end{center}\end{minipage}};

            \draw (15,60) node  [font=\normalsize] [align=left] {\begin{minipage}[lt]{20.4pt}\setlength\topsep{0pt}
            \begin{center}
            $\displaystyle \x$
            \end{center}\end{minipage}};

            \draw (15,160) node  [font=\normalsize] [align=left] {\begin{minipage}[lt]{20.4pt}\setlength\topsep{0pt}
            \begin{center}
            $\displaystyle \z$
            \end{center}\end{minipage}};
        \end{tikzpicture}
        \caption{ 
            Vertices are divided into two groups:
            the upper vertices  correspond to $\x$
            while the lower ones correspond to $\z$.
            }
        \label{fig:example-transformed-sparsity-before}
    \end{subfigure}
    \hfill
    \begin{subfigure}[t]{0.45\textwidth}
        \centering
        \tikzset{every picture/.style={line width=0.75pt}}
        \begin{tikzpicture}[x=0.75pt,y=0.75pt,yscale=-0.6,xscale=0.6]
            \draw   (215,110) -- (315,70) -- (388,125.35) -- (141,53.35) -- cycle ;
            \draw  [dash pattern={on 4.5pt off 4.5pt}]  (115,40) -- (315,140) ;
            \draw  [dash pattern={on 4.5pt off 4.5pt}]  (115,140) -- (315,40) ;
            \draw  [dash pattern={on 0.84pt off 2.51pt}]  (315,40) -- (315,140) ;
            \draw  [dash pattern={on 0.84pt off 2.51pt}]  (115,40) -- (115,140) ;
            \draw  [dash pattern={on 0.84pt off 2.51pt}]  (215,40) -- (215,140) ;
            \draw  [dash pattern={on 0.84pt off 2.51pt}]  (415,40) -- (415,140) ;
            \draw  [fill={rgb, 255:red, 255; green, 255; blue, 255 }  ,fill opacity=1 ] (85,40) .. controls (85,23.43) and (98.43,10) .. (115,10) .. controls (131.57,10) and (145,23.43) .. (145,40) .. controls (145,56.57) and (131.57,70) .. (115,70) .. controls (98.43,70) and (85,56.57) .. (85,40) -- cycle ;

            \draw  [fill={rgb, 255:red, 255; green, 255; blue, 255 }  ,fill opacity=1 ] (285,40) .. controls (285,23.43) and (298.43,10) .. (315,10) .. controls (331.57,10) and (345,23.43) .. (345,40) .. controls (345,56.57) and (331.57,70) .. (315,70) .. controls (298.43,70) and (285,56.57) .. (285,40) -- cycle ;
            \draw  [fill={rgb, 255:red, 255; green, 255; blue, 255 }  ,fill opacity=1 ] (185,40) .. controls (185,23.43) and (198.43,10) .. (215,10) .. controls (231.57,10) and (245,23.43) .. (245,40) .. controls (245,56.57) and (231.57,70) .. (215,70) .. controls (198.43,70) and (185,56.57) .. (185,40) -- cycle ;
            \draw  [fill={rgb, 255:red, 255; green, 255; blue, 255 }  ,fill opacity=1 ] (385,40) .. controls (385,23.43) and (398.43,10) .. (415,10) .. controls (431.57,10) and (445,23.43) .. (445,40) .. controls (445,56.57) and (431.57,70) .. (415,70) .. controls (398.43,70) and (385,56.57) .. (385,40) -- cycle ;
            \draw  [fill={rgb, 255:red, 255; green, 255; blue, 255 }  ,fill opacity=1 ] (85,140) .. controls (85,123.43) and (98.43,110) .. (115,110) .. controls (131.57,110) and (145,123.43) .. (145,140) .. controls (145,156.57) and (131.57,170) .. (115,170) .. controls (98.43,170) and (85,156.57) .. (85,140) -- cycle ;
            \draw  [fill={rgb, 255:red, 255; green, 255; blue, 255 }  ,fill opacity=1 ] (285,140) .. controls (285,123.43) and (298.43,110) .. (315,110) .. controls (331.57,110) and (345,123.43) .. (345,140) .. controls (345,156.57) and (331.57,170) .. (315,170) .. controls (298.43,170) and (285,156.57) .. (285,140) -- cycle ;
            \draw  [fill={rgb, 255:red, 255; green, 255; blue, 255 }  ,fill opacity=1 ] (185,140) .. controls (185,123.43) and (198.43,110) .. (215,110) .. controls (231.57,110) and (245,123.43) .. (245,140) .. controls (245,156.57) and (231.57,170) .. (215,170) .. controls (198.43,170) and (185,156.57) .. (185,140) -- cycle ;
            \draw  [fill={rgb, 255:red, 255; green, 255; blue, 255 }  ,fill opacity=1 ] (385,140) .. controls (385,123.43) and (398.43,110) .. (415,110) .. controls (431.57,110) and (445,123.43) .. (445,140) .. controls (445,156.57) and (431.57,170) .. (415,170) .. controls (398.43,170) and (385,156.57) .. (385,140) -- cycle ;

            \draw (215,40) node  [font=\large] [align=left] {\begin{minipage}[lt]{40.8pt}\setlength\topsep{0pt}
            \begin{center}
            {\fontfamily{pcr}\selectfont 6}
            \end{center}\end{minipage}};

            \draw (415,40) node  [font=\large] [align=left] {\begin{minipage}[lt]{40.8pt}\setlength\topsep{0pt}
            \begin{center}
            {\fontfamily{pcr}\selectfont 8}
            \end{center}\end{minipage}};

            \draw (115,40) node  [font=\large] [align=left] {\begin{minipage}[lt]{40.8pt}\setlength\topsep{0pt}
            \begin{center}
            {\fontfamily{pcr}\selectfont 1}
            \end{center}\end{minipage}};

            \draw (315,40) node  [font=\large] [align=left] {\begin{minipage}[lt]{40.8pt}\setlength\topsep{0pt}
            \begin{center}
            {\fontfamily{pcr}\selectfont 3}
            \end{center}\end{minipage}};

            \draw (215,140) node  [font=\large] [align=left] {\begin{minipage}[lt]{40.8pt}\setlength\topsep{0pt}
            \begin{center}
            {\fontfamily{pcr}\selectfont 2}
            \end{center}\end{minipage}};

            \draw (415,140) node  [font=\large] [align=left] {\begin{minipage}[lt]{40.8pt}\setlength\topsep{0pt}
            \begin{center}
            {\fontfamily{pcr}\selectfont 4}
            \end{center}\end{minipage}};

            \draw (115,140) node  [font=\large] [align=left] {\begin{minipage}[lt]{40.8pt}\setlength\topsep{0pt}
            \begin{center}
            {\fontfamily{pcr}\selectfont 5}
            \end{center}\end{minipage}};

            \draw (315,140) node  [font=\large] [align=left] {\begin{minipage}[lt]{40.8pt}\setlength\topsep{0pt}
            \begin{center}
            {\fontfamily{pcr}\selectfont 7}
            \end{center}\end{minipage}};

            \draw (30,40) node  [font=\normalsize] [align=left] {\begin{minipage}[lt]{54.4pt}\setlength\topsep{0pt}
            \begin{center}
            $\displaystyle L_{1} \cup L_{2}$
            \end{center}\end{minipage}};

            \draw (30,140) node  [font=\normalsize] [align=left] {\begin{minipage}[lt]{54.4pt}\setlength\topsep{0pt}
            \begin{center}
            $\displaystyle R_{1} \cup R_{2}$
            \end{center}\end{minipage}};
        \end{tikzpicture}
        \caption{
            Vertices are reorganized to show the bipartite structure of the graph.
            }
        \label{fig:example-transformed-sparsity-after}
    \end{subfigure}
    \caption{
        Aggregated sparsity pattern graph of the transformed example.
        The solid lines and the dashed lines come from $N^p_+$ and $N^p_-$, respectively.
        The dotted lines are for the new constraint $\|\x + \z\|^2 \leq 0$.
        }
    \label{fig:example-transformed-sparsity}
\end{figure}

\begin{example}
Now, consider an instance of QCQP~\eqref{eq:hqcqp}
with $n=4$, $Q^p_{24} = 0 \, (p \in [0, m])$ and the edge signs as:
\begin{equation*}
    \begin{aligned}
        \sigma_{12} &= +1, & \sigma_{13} &= -1, & \sigma_{14} &= +1, & \sigma_{23} &= +1, & \sigma_{34} &= +1.
    \end{aligned}
\end{equation*}
\autoref{fig:example-edge-signs} illustrates the above signs.
We also suppose that $Q^0_{ij} \neq 0$ for all $(i, j) \in \EC$.
Then, for any distinct $i, j \in [n]$,
the set $\left\{Q^0_{ij}, \ldots, Q^m_{ij}\right\}$ is sign-definite by definition.
Since there exist odd cycles, e.g., $\{(1, 2), (2, 3), (3, 1)\}$,
the aggregated sparsity pattern graph of a QCQP with the above edge signs is not bipartite.
Next, we transform the QCQP instance into an equivalent QCQP with bipartite structures.
Since $n=4$, we see $\bar{\VC} = [8]$.
\autoref{fig:example-transformed-sparsity-before} displays $\bar{G}$ from
\begin{equation*}
    \begin{bmatrix} D^p + 2N^p_+ & N^p_- \\ N^p_- & O \end{bmatrix} =
    \left[
        \begin{array}{c|c}
            \begin{matrix}
                Q^p_{11} & Q^p_{12} & 0 & Q^p_{14} \\
                Q^p_{21} & Q^p_{22} & Q^p_{23} & 0 \\
                0 & Q^p_{32} & Q^p_{33} & Q^p_{34} \\
                Q^p_{41} & 0 & Q^p_{43} & Q^p_{44}
            \end{matrix} &
            \begin{matrix}
                0 & 0 & -\frac{1}{2}Q^p_{13} & 0 \\
                0 & 0 & 0 & 0 \\
                -\frac{1}{2}Q^p_{31} & 0 & 0 & 0 \\
                0 & 0 & 0 & 0
            \end{matrix} \\ \hline
            \begin{matrix}
                0 & 0 & -\frac{1}{2}Q^p_{13} & 0 \\
                0 & 0 & 0 & 0 \\
                -\frac{1}{2}Q^p_{31} & 0 & 0 & 0 \\
                0 & 0 & 0 & 0
            \end{matrix} & O \\
        \end{array}
    \right] +
    \delta \left[
        \begin{array}{c|c}
            2 I & I \\ \hline
            I & O
        \end{array}
    \right]
\end{equation*}
and $[I\; I; I\; I]$.
There exist three types of edges:
\begin{equation*}
    \def\arraystretch{1.5}
    \left\{\begin{array}{rl}
        \text{  (i)} & (1, 2), (2, 3), (3, 4), (1, 4); \\
        \text{ (ii)} & (1, 7), (3, 5); \\
        \text{(iii)} & (1, 5), (2, 6), (3, 7), (4, 8).
    \end{array}\right.
\end{equation*}
The edges in (i) and (ii) are derived from
four $N^p_+$ on the upper-left of the data matrices,
and two $N^p_-$ on the upper-right and the lower-left of the data matrices, respectively.
The edges for (iii) represents off-diagonal elements in $[I\; I; I\; I]$ in the new constraint.
In \autoref{fig:example-transformed-sparsity-before},
the cycle in the solid lines is bipartite with the vertices $\{1,2,3,4\}$,
and hence its vertices can be divided into two distinct sets $L_1 = \{1, 3\}$ and $R_1 = \{2, 4\}$.
If we let $L_2 \coloneqq \{6, 8\}$ and $R_2 \coloneqq \{5, 7\}$,
there are no edges between any distinct $i, j$ in $L_1 \cup L_2$,
and the same is true for $R_1 \cup R_2$.
The graph $\bar{G}$ is thus bipartite (\autoref{fig:example-transformed-sparsity-after}).
We can conclude that the SDP relaxation of \eqref{eq:decomposed-hqcqp} is exact
by Corollary~\ref{coro:nonnegative-offdiagonal-connected}.
\end{example}

Similarly, the SDP relaxation of any QCQP that satisfies Theorem~\ref{thm:sojoudi-theorem}
can be shown to be exact by the transformation.
Therefore,
Theorem~\ref{thm:system-based-condition-connected} includes a wider classes of QCQPs than Theorem~\ref{thm:sojoudi-theorem}.
We prove this assertion in the following.
%
\begin{prop}
    \label{prop:weaker-than-sojoudi-connected}
    Suppose that Assumption~\ref{assum:new-assumption} holds,
    the aggregated sparsity pattern graph $G(\VC, \EC)$ of \eqref{eq:hqcqp}
    is connected, 
    and for all $(i, j) \in \EC$, $Q^0_{ij} \neq 0$.
    If \eqref{eq:hqcqp} satisfies the assumption of Theorem~\ref{thm:sojoudi-theorem},
    then \eqref{eq:hqcqp} also satisfies that of Corollary~\ref{coro:nonnegative-offdiagonal-connected}.
    In addition, the exactness of its SDP relaxation
    can be proved by Theorem~\ref{thm:system-based-condition-connected}.
\end{prop}
\begin{proof}
    Let $\bar{G}(\bar{\VC}, \bar{\EC})$ be
    the aggregated sparsity pattern graph of \eqref{eq:decomposed-hqcqp}.
    Since the number of variables is $2n$, $\bar{\VC} = [2n]$ holds.
    The edges in $\bar{G}$ are:
    \begin{equation*}
        \def\arraystretch{1.5}
        \left\{\begin{array}{rll}
            \text{  (i)} & (i, j) & \text{for $i, j \in \VC$ such that $\sigma_{ij} = +1$}, \\
            \text{ (ii)} & (i, j + n), (j, i + n) &
                \text{for $i, j \in \VC$ such that $\sigma_{ij} = -1$}, \\
            \text{(iii)} & (i, i + n) & \text{for $i \in \VC$}.
        \end{array}\right.
    \end{equation*}
    Note that
    no edges exist among the vertices in $\{n+1, \ldots, 2n\}$.
    By the definition of \eqref{eq:decomposed-hqcqp},
    an edge $(i, j)$ with  $\sigma_{ij} = -1$ in $G$
    is decomposed into two paths with positive signs in $\bar{G}$:
    (a) the edges $(j, i + n)$ and $(i + n, i)$;
    (b) the edges $(i, j + n)$ and $(j + n, j)$, as shown in \autoref{fig:transform-minus-edge-sign}.
    Since $G$ is connected, so is the graph $\bar{G}$. 
    Recall that all off-diagonal elements of
    the data matrices in \eqref{eq:decomposed-hqcqp} are nonnegative, since both $N^p_+$ and $N^p_-$ are nonnegative matrices.
    In particular, for each $(i, j) \in \bar{\EC}$,
    the $(i, j)$th element of the matrix in the objective function is not only nonnegative but also positive by assumption.
    Thus, to apply Corollary~\ref{coro:nonnegative-offdiagonal-connected}, it remains to show that $\bar{G}$ is bipartite.


    Assume on the contrary there exists an odd cycle $\bar{\CC}$ in $\bar{G}$.
    Let $\bar{\UC} \subseteq [n+1, 2n]$ denote the set of vertices on $[n+1, 2n]$ in $\bar{\CC}$.
    As illustrated in \autoref{fig:transform-minus-edge-sign},
    any vertex $v \coloneqq i + n \in \bar{\UC}$ connects with $i$ and $j \in \VC$ in $\bar{\CC}$.
    Hence for every vertex $v \in \bar{\UC}$,
    by removing the edges $(i, v)$ and $(v, j)$ from $\bar{\CC}$ and
    adding the edge $(i, j)$ with the negative sign to $\bar{\CC}$,
    we obtain a new cycle $\CC$ in $G$.
    Since $2|\bar{\UC}|$ edges are removed and $|\bar{\UC}|$ edges are added in this procedure,
    it follows $|\CC| = |\bar{\CC}| - 2|\bar{\UC}| + |\bar{\UC}| = |\bar{\CC}| - |\bar{\UC}|$.
    \autoref{fig:removing-adding-cycle-edges} displays a case for $|\bar{\UC}| = 2$.
    Thus, if $|\bar{\UC}|$ is even (odd), $|\CC|$ is odd (resp., even),
    hence, by \eqref{eq:sign-constraint-simple-cycle} in Theorem~\ref{thm:sojoudi-theorem},
    the number of negative edges in $\CC$ must be odd (resp., even).
    However, the number of negative edges in $\CC$ is equal to $|\bar{\UC}|$
    since $\bar{\CC}$ has no negative edges and all the additional edges
    in the conversion from $\bar{\CC}$ to $\CC$
    are negative.
    This is a contradiction.
    Therefore, there are no odd cycles in $\bar{G}$,
    which implies $\bar{G}$ is bipartite.
    Since \eqref{eq:decomposed-hqcqp} satisfies the assumptions of Corollary~\ref{coro:nonnegative-offdiagonal-connected},
    it also satisfies the assumptions of Theorem~\ref{thm:system-based-condition-connected}.
\end{proof}

\begin{figure}[t]
    \centering
    \tikzset{every picture/.style={line width=0.75pt}} 
    \begin{minipage}[t]{0.54\textwidth}
        \centering
        \begin{tikzpicture}[x=0.75pt,y=0.75pt,yscale=-0.6,xscale=0.6]
            \draw  [draw opacity=0] (324.92,109.3) .. controls (318.07,127.62) and (305,140) .. (290,140) .. controls (267.91,140) and (250,113.14) .. (250,80) .. controls (250,46.86) and (267.91,20) .. (290,20) .. controls (305,20) and (318.07,32.38) .. (324.92,50.7) -- (290,80) -- cycle ; \draw   (324.92,109.3) .. controls (318.07,127.62) and (305,140) .. (290,140) .. controls (267.91,140) and (250,113.14) .. (250,80) .. controls (250,46.86) and (267.91,20) .. (290,20) .. controls (305,20) and (318.07,32.38) .. (324.92,50.7) ;
            \draw [color={rgb, 255:red, 0; green, 0; blue, 0 }  ,draw opacity=1 ][line width=1.5]  [dash pattern={on 1.69pt off 2.76pt}]  (324.92,50.7) -- (404.9,50.7) ;
            \draw  [draw opacity=0][fill={rgb, 255:red, 0; green, 0; blue, 0 }  ,fill opacity=1 ] (318.92,50.7) .. controls (318.92,47.39) and (321.6,44.7) .. (324.92,44.7) .. controls (328.23,44.7) and (330.92,47.39) .. (330.92,50.7) .. controls (330.92,54.02) and (328.23,56.7) .. (324.92,56.7) .. controls (321.6,56.7) and (318.92,54.02) .. (318.92,50.7) -- cycle ;
            \draw  [draw opacity=0][fill={rgb, 255:red, 0; green, 0; blue, 0 }  ,fill opacity=1 ] (318.92,109.3) .. controls (318.92,105.98) and (321.6,103.3) .. (324.92,103.3) .. controls (328.23,103.3) and (330.92,105.98) .. (330.92,109.3) .. controls (330.92,112.61) and (328.23,115.3) .. (324.92,115.3) .. controls (321.6,115.3) and (318.92,112.61) .. (318.92,109.3) -- cycle ;
            \draw [color={rgb, 255:red, 0; green, 0; blue, 0 }  ,draw opacity=1 ][line width=1.5]  [dash pattern={on 1.69pt off 2.76pt}]  (324.92,109.3) -- (404.9,109.3) ;
            \draw [color={rgb, 255:red, 0; green, 0; blue, 0 }  ,draw opacity=1 ][line width=1.5]  [dash pattern={on 5.63pt off 4.5pt}]  (324.92,50.7) -- (404.9,109.3) ;
            \draw [color={rgb, 255:red, 0; green, 0; blue, 0 }  ,draw opacity=1 ][line width=1.5]  [dash pattern={on 5.63pt off 4.5pt}]  (324.92,109.3) -- (404.9,50.7) ;
            \draw  [draw opacity=0][fill={rgb, 255:red, 0; green, 0; blue, 0 }  ,fill opacity=1 ] (398.9,50.7) .. controls (398.9,47.39) and (401.59,44.7) .. (404.9,44.7) .. controls (408.21,44.7) and (410.9,47.39) .. (410.9,50.7) .. controls (410.9,54.01) and (408.21,56.7) .. (404.9,56.7) .. controls (401.59,56.7) and (398.9,54.01) .. (398.9,50.7) -- cycle ;
            \draw  [draw opacity=0][fill={rgb, 255:red, 0; green, 0; blue, 0 }  ,fill opacity=1 ] (398.9,109.3) .. controls (398.9,105.98) and (401.59,103.3) .. (404.9,103.3) .. controls (408.21,103.3) and (410.9,105.98) .. (410.9,109.3) .. controls (410.9,112.61) and (408.21,115.3) .. (404.9,115.3) .. controls (401.59,115.3) and (398.9,112.61) .. (398.9,109.3) -- cycle ;

            \draw  [draw opacity=0] (114.92,109.3) .. controls (108.07,127.62) and (95,140) .. (80,140) .. controls (57.91,140) and (40,113.14) .. (40,80) .. controls (40,46.86) and (57.91,20) .. (80,20) .. controls (95,20) and (108.07,32.38) .. (114.92,50.7) -- (80,80) -- cycle ; \draw   (114.92,109.3) .. controls (108.07,127.62) and (95,140) .. (80,140) .. controls (57.91,140) and (40,113.14) .. (40,80) .. controls (40,46.86) and (57.91,20) .. (80,20) .. controls (95,20) and (108.07,32.38) .. (114.92,50.7) ;
            \draw  [draw opacity=0][fill={rgb, 255:red, 0; green, 0; blue, 0 }  ,fill opacity=1 ] (108.92,50.7) .. controls (108.92,47.39) and (111.6,44.7) .. (114.92,44.7) .. controls (118.23,44.7) and (120.92,47.39) .. (120.92,50.7) .. controls (120.92,54.02) and (118.23,56.7) .. (114.92,56.7) .. controls (111.6,56.7) and (108.92,54.02) .. (108.92,50.7) -- cycle ;
            \draw  [draw opacity=0][fill={rgb, 255:red, 0; green, 0; blue, 0 }  ,fill opacity=1 ] (108.92,109.3) .. controls (108.92,105.98) and (111.6,103.3) .. (114.92,103.3) .. controls (118.23,103.3) and (120.92,105.98) .. (120.92,109.3) .. controls (120.92,112.61) and (118.23,115.3) .. (114.92,115.3) .. controls (111.6,115.3) and (108.92,112.61) .. (108.92,109.3) -- cycle ;
            \draw [color={rgb, 255:red, 0; green, 0; blue, 0 }  ,draw opacity=1 ][line width=1.5]    (114.92,50.7) -- (114.92,109.3) ;
            \draw   (150,80) -- (167.5,60) -- (167.5,70) -- (202.5,70) -- (202.5,60) -- (220,80) -- (202.5,100) -- (202.5,90) -- (167.5,90) -- (167.5,100) -- cycle ;

            \draw (290,157.5) node  [font=\normalsize] [align=left] {\begin{minipage}[lt]{68pt}\setlength\topsep{0pt}
            \begin{center}
            $\displaystyle \mathcal{C} \setminus \{( i,j)\}$
            \end{center}\end{minipage}};
            \draw (240,30) node  [font=\large] [align=left] {\begin{minipage}[lt]{27.2pt}\setlength\topsep{0pt}
            \begin{center}
            $\displaystyle \overline{G}$
            \end{center}\end{minipage}};
            \draw (404.9,34.7) node  [font=\normalsize] [align=left] {\begin{minipage}[lt]{40.8pt}\setlength\topsep{0pt}
            \begin{center}
            $\displaystyle i+n$
            \end{center}\end{minipage}};
            \draw (404.9,129.3) node  [font=\normalsize] [align=left] {\begin{minipage}[lt]{40.8pt}\setlength\topsep{0pt}
            \begin{center}
            $\displaystyle j+n$
            \end{center}\end{minipage}};
            \draw (334.92,34.7) node  [font=\normalsize] [align=left] {\begin{minipage}[lt]{13.6pt}\setlength\topsep{0pt}
            \begin{center}
            $\displaystyle i$
            \end{center}\end{minipage}};
            \draw (334.92,129.3) node  [font=\normalsize] [align=left] {\begin{minipage}[lt]{13.6pt}\setlength\topsep{0pt}
            \begin{center}
            $\displaystyle j$
            \end{center}\end{minipage}};
            \draw (80,157.5) node  [font=\normalsize] [align=left] {\begin{minipage}[lt]{68pt}\setlength\topsep{0pt}
            \begin{center}
            $\displaystyle \mathcal{C}$
            \end{center}\end{minipage}};
            \draw (30,30) node  [font=\large] [align=left] {\begin{minipage}[lt]{27.2pt}\setlength\topsep{0pt}
            \begin{center}
            $\displaystyle G$
            \end{center}\end{minipage}};
            \draw (124.92,34.7) node  [font=\normalsize] [align=left] {\begin{minipage}[lt]{13.6pt}\setlength\topsep{0pt}
            \begin{center}
            $\displaystyle i$
            \end{center}\end{minipage}};
            \draw (124.92,129.3) node  [font=\normalsize] [align=left] {\begin{minipage}[lt]{13.6pt}\setlength\topsep{0pt}
            \begin{center}
            $\displaystyle j$
            \end{center}\end{minipage}};
            \draw (100,80) node  [font=\large] [align=left] {\begin{minipage}[lt]{20.4pt}\setlength\topsep{0pt}
            \begin{center}
            $\displaystyle -$
            \end{center}\end{minipage}};
        \end{tikzpicture}
        \caption{
            An edge with the negative sign.
            If the cycle $\CC$ has the edge $(i, j)$ with $\sigma_{ij} = -1$,
            then $(i, j)$ is decomposed into two paths:
            (a) $(j, i+n)$ and $(i+n, i)$ via the vertex $i+n$;
            (b) $(i, j+n)$ and $(j+n, j)$ via the vertex $j+n$.
            }
        \label{fig:transform-minus-edge-sign}
    \end{minipage}
    \hfill
    \begin{minipage}[t]{0.42\textwidth}
        \centering
        \begin{tikzpicture}[x=0.75pt,y=0.75pt,yscale=-0.725,xscale=0.725]
            \draw  [draw opacity=0][fill={rgb, 255:red, 0; green, 0; blue, 0 }  ,fill opacity=1 ] (45,10) .. controls (45,7.24) and (47.24,5) .. (50,5) .. controls (52.76,5) and (55,7.24) .. (55,10) .. controls (55,12.76) and (52.76,15) .. (50,15) .. controls (47.24,15) and (45,12.76) .. (45,10) -- cycle ;
            \draw  [draw opacity=0][fill={rgb, 255:red, 0; green, 0; blue, 0 }  ,fill opacity=1 ] (45,110) .. controls (45,107.24) and (47.24,105) .. (50,105) .. controls (52.76,105) and (55,107.24) .. (55,110) .. controls (55,112.76) and (52.76,115) .. (50,115) .. controls (47.24,115) and (45,112.76) .. (45,110) -- cycle ;
            \draw  [draw opacity=0][fill={rgb, 255:red, 0; green, 0; blue, 0 }  ,fill opacity=1 ] (45,40) .. controls (45,37.24) and (47.24,35) .. (50,35) .. controls (52.76,35) and (55,37.24) .. (55,40) .. controls (55,42.76) and (52.76,45) .. (50,45) .. controls (47.24,45) and (45,42.76) .. (45,40) -- cycle ;
            \draw  [draw opacity=0][fill={rgb, 255:red, 0; green, 0; blue, 0 }  ,fill opacity=1 ] (45,80) .. controls (45,77.24) and (47.24,75) .. (50,75) .. controls (52.76,75) and (55,77.24) .. (55,80) .. controls (55,82.76) and (52.76,85) .. (50,85) .. controls (47.24,85) and (45,82.76) .. (45,80) -- cycle ;
            \draw    (10,10) -- (50,10) -- (90,25) -- (50,40) -- (50,80) -- (90,95) -- (50,110) -- (10,110) -- (10,10);

            \draw  [draw opacity=0][fill={rgb, 255:red, 255; green, 255; blue, 255 }  ,fill opacity=1 ] (45,64) -- (55,59) -- (55,55) -- (45,60) -- cycle ;
            \draw    (45,60) -- (55,55) ;
            \draw    (45,64) -- (55,59) ;
            \draw  [draw opacity=0][fill={rgb, 255:red, 255; green, 255; blue, 255 }  ,fill opacity=1 ] (5,64) -- (15,59) -- (15,55) -- (5,60) -- cycle ;
            \draw    (5,60) -- (15,55) ;
            \draw    (5,64) -- (15,59) ;

            \draw  [color={rgb, 255:red, 0; green, 0; blue, 0 }  ,draw opacity=1 ][fill={rgb, 255:red, 255; green, 255; blue, 255 }  ,fill opacity=1 ] (85,25) .. controls (85,22.24) and (87.24,20) .. (90,20) .. controls (92.76,20) and (95,22.24) .. (95,25) .. controls (95,27.76) and (92.76,30) .. (90,30) .. controls (87.24,30) and (85,27.76) .. (85,25) -- cycle ;
            \draw  [color={rgb, 255:red, 0; green, 0; blue, 0 }  ,draw opacity=1 ][fill={rgb, 255:red, 255; green, 255; blue, 255 }  ,fill opacity=1 ] (85,95) .. controls (85,92.24) and (87.24,90) .. (90,90) .. controls (92.76,90) and (95,92.24) .. (95,95) .. controls (95,97.76) and (92.76,100) .. (90,100) .. controls (87.24,100) and (85,97.76) .. (85,95) -- cycle ;
            \draw  [draw opacity=0][fill={rgb, 255:red, 0; green, 0; blue, 0 }  ,fill opacity=1 ] (125,10) .. controls (125,7.24) and (127.24,5) .. (130,5) .. controls (132.76,5) and (135,7.24) .. (135,10) .. controls (135,12.76) and (132.76,15) .. (130,15) .. controls (127.24,15) and (125,12.76) .. (125,10) -- cycle ;
            \draw  [draw opacity=0][fill={rgb, 255:red, 0; green, 0; blue, 0 }  ,fill opacity=1 ] (125,110) .. controls (125,107.24) and (127.24,105) .. (130,105) .. controls (132.76,105) and (135,107.24) .. (135,110) .. controls (135,112.76) and (132.76,115) .. (130,115) .. controls (127.24,115) and (125,112.76) .. (125,110) -- cycle ;
            \draw  [draw opacity=0][fill={rgb, 255:red, 0; green, 0; blue, 0 }  ,fill opacity=1 ] (125,40) .. controls (125,37.24) and (127.24,35) .. (130,35) .. controls (132.76,35) and (135,37.24) .. (135,40) .. controls (135,42.76) and (132.76,45) .. (130,45) .. controls (127.24,45) and (125,42.76) .. (125,40) -- cycle ;
            \draw  [draw opacity=0][fill={rgb, 255:red, 0; green, 0; blue, 0 }  ,fill opacity=1 ] (125,80) .. controls (125,77.24) and (127.24,75) .. (130,75) .. controls (132.76,75) and (135,77.24) .. (135,80) .. controls (135,82.76) and (132.76,85) .. (130,85) .. controls (127.24,85) and (125,82.76) .. (125,80) -- cycle ;
            \draw    (130,10) -- (170,25) -- (130,40) ;
            \draw    (130,80) -- (170,95) -- (130,110) ;
            \draw  [color={rgb, 255:red, 0; green, 0; blue, 0 }  ,draw opacity=1 ][fill={rgb, 255:red, 255; green, 255; blue, 255 }  ,fill opacity=1 ] (165,25) .. controls (165,22.24) and (167.24,20) .. (170,20) .. controls (172.76,20) and (175,22.24) .. (175,25) .. controls (175,27.76) and (172.76,30) .. (170,30) .. controls (167.24,30) and (165,27.76) .. (165,25) -- cycle ;
            \draw  [color={rgb, 255:red, 0; green, 0; blue, 0 }  ,draw opacity=1 ][fill={rgb, 255:red, 255; green, 255; blue, 255 }  ,fill opacity=1 ] (165,95) .. controls (165,92.24) and (167.24,90) .. (170,90) .. controls (172.76,90) and (175,92.24) .. (175,95) .. controls (175,97.76) and (172.76,100) .. (170,100) .. controls (167.24,100) and (165,97.76) .. (165,95) -- cycle ;

            \draw  [draw opacity=0][fill={rgb, 255:red, 0; green, 0; blue, 0 }  ,fill opacity=1 ] (215,10) .. controls (215,7.24) and (217.24,5) .. (220,5) .. controls (222.76,5) and (225,7.24) .. (225,10) .. controls (225,12.76) and (222.76,15) .. (220,15) .. controls (217.24,15) and (215,12.76) .. (215,10) -- cycle ;
            \draw  [draw opacity=0][fill={rgb, 255:red, 0; green, 0; blue, 0 }  ,fill opacity=1 ] (215,110) .. controls (215,107.24) and (217.24,105) .. (220,105) .. controls (222.76,105) and (225,107.24) .. (225,110) .. controls (225,112.76) and (222.76,115) .. (220,115) .. controls (217.24,115) and (215,112.76) .. (215,110) -- cycle ;
            \draw  [draw opacity=0][fill={rgb, 255:red, 0; green, 0; blue, 0 }  ,fill opacity=1 ] (215,40) .. controls (215,37.24) and (217.24,35) .. (220,35) .. controls (222.76,35) and (225,37.24) .. (225,40) .. controls (225,42.76) and (222.76,45) .. (220,45) .. controls (217.24,45) and (215,42.76) .. (215,40) -- cycle ;
            \draw  [draw opacity=0][fill={rgb, 255:red, 0; green, 0; blue, 0 }  ,fill opacity=1 ] (215,80) .. controls (215,77.24) and (217.24,75) .. (220,75) .. controls (222.76,75) and (225,77.24) .. (225,80) .. controls (225,82.76) and (222.76,85) .. (220,85) .. controls (217.24,85) and (215,82.76) .. (215,80) -- cycle ;
            \draw    (220,10) -- (220,40) ;
            \draw    (220,80) -- (220,110) ;

            \draw  [draw opacity=0][fill={rgb, 255:red, 0; green, 0; blue, 0 }  ,fill opacity=1 ] (320,10) .. controls (320,7.24) and (322.24,5) .. (325,5) .. controls (327.76,5) and (330,7.24) .. (330,10) .. controls (330,12.76) and (327.76,15) .. (325,15) .. controls (322.24,15) and (320,12.76) .. (320,10) -- cycle ;
            \draw  [draw opacity=0][fill={rgb, 255:red, 0; green, 0; blue, 0 }  ,fill opacity=1 ] (320,110) .. controls (320,107.24) and (322.24,105) .. (325,105) .. controls (327.76,105) and (330,107.24) .. (330,110) .. controls (330,112.76) and (327.76,115) .. (325,115) .. controls (322.24,115) and (320,112.76) .. (320,110) -- cycle ;
            \draw  [draw opacity=0][fill={rgb, 255:red, 0; green, 0; blue, 0 }  ,fill opacity=1 ] (320,40) .. controls (320,37.24) and (322.24,35) .. (325,35) .. controls (327.76,35) and (330,37.24) .. (330,40) .. controls (330,42.76) and (327.76,45) .. (325,45) .. controls (322.24,45) and (320,42.76) .. (320,40) -- cycle ;
            \draw  [draw opacity=0][fill={rgb, 255:red, 0; green, 0; blue, 0 }  ,fill opacity=1 ] (320,80) .. controls (320,77.24) and (322.24,75) .. (325,75) .. controls (327.76,75) and (330,77.24) .. (330,80) .. controls (330,82.76) and (327.76,85) .. (325,85) .. controls (322.24,85) and (320,82.76) .. (320,80) -- cycle ;
            \draw    (285,10) -- (285,110) -- (325,110) -- (325,10) -- (285,10) ;

            \draw  [draw opacity=0][fill={rgb, 255:red, 255; green, 255; blue, 255 }  ,fill opacity=1 ] (320,64) -- (330,59) -- (330,55) -- (320,60) -- cycle ;
            \draw    (320,60) -- (330,55) ;
            \draw    (320,64) -- (330,59) ;

            \draw  [draw opacity=0][fill={rgb, 255:red, 255; green, 255; blue, 255 }  ,fill opacity=1 ] (280,64) -- (290,59) -- (290,55) -- (280,60) -- cycle ;
            \draw    (280,60) -- (290,55) ;
            \draw    (280,64) -- (290,59) ;

            \draw (110,60) node   [align=left] {$\displaystyle -$};
            \draw (195,60) node   [align=left] {$\displaystyle +$};
            \draw (250,60) node   [align=left] {$\displaystyle =$};
            \draw (50,135) node   [align=left] {$\displaystyle | \overline{\mathcal{C}}| $};
            \draw (150,135) node   [align=left] {$\displaystyle 2| \overline{U}| $};
            \draw (220,135) node   [align=left] {$\displaystyle | \overline{U}| $};
            \draw (305,135) node   [align=left] {$\displaystyle | \mathcal{C}| $};
        \end{tikzpicture}
        \caption{
            Removing and adding edges, and calculating of the number of edges if $\bar{\UC} = 2$.
            The black circles are the vertices in $[n]$ while the white circles represent those in $[n+1, 2n]$.
            }
        \label{fig:removing-adding-cycle-edges}
    \end{minipage}
\end{figure}

\noindent
Proposition~\ref{prop:weaker-than-sojoudi-connected} is proved 
under the assumptions that: (i) $G$ is connected; (ii) for all $(i, j) \in \EC$, $Q^0_{ij} \neq 0$.
These assumptions may seem strong;
however, we will show that they can be removed
 using 
Corollary~\ref{coro:nonnegative-offdiagonal} in section 4. 


At the end of this section, we apply Proposition~\ref{prop:weaker-than-sojoudi-connected} to
a class of QCQPs where all the off-diagonal elements of every matrix $Q^0, \ldots, Q^m$ are nonpositive.
We call QCQPs in this class nonpositive off-diagonal QCQPs.
It is well-known that their SDP relaxations are exact~\cite{kim2003exact}.
By applying the same transformation above,
 we obtain \eqref{eq:decomposed-hqcqp} with $N^p_+ = O$ for every $p$
since no positive off-diagonal elements exist.
The diagonal elements of $D^p$ do not generate edges in the aggregated sparsity pattern graph, thus,
the 
data matrices in \eqref{eq:decomposed-hqcqp} induce a bipartite sparsity pattern graph.
Therefore, the SDP relaxation is exact.
This can be regarded as an alternative proof for \cite{kim2003exact} and Corollary~\ref{coro:sojoudi-corollary1}\ref{cond:sojoudi-arbitrary}.
\begin{coro} \label{coro:nonpositive-offdiagonal-connected}
    Under Assumption~\ref{assum:new-assumption},
    the SDP relaxation of a nonpositive off-diagonal QCQP is exact
    if the aggregate spartiy pattern graph $G(\VC, \EC)$ of \eqref{eq:hqcqp} is connected
    and $Q^0_{ij} < 0$ for all $(i, j) \in \EC$.
\end{coro}

\section{Perturbation for disconnected aggregated sparsity pattern graph} \label{sec:perturbation}

The connectivity of $G$ has played an important role  for
our main theorem in \autoref{sec:main}.
For QCQPs with sparse data matrices,
the connectivity assumption might be a difficult condition to be satisfied.
In this section,
we replace the assumption for connected graphs
by a slightly different assumption (Assumption~\ref{assum:new-assumption-strong}),
and present  a new condition for the exact SDP relaxation.

The following assumption is slightly stronger than Assumption~\ref{assum:new-assumption}
in the sense that it requires the existence of a feasible interior point of \eqref{eq:hsdrd}.
However, it can be satisfied in practice without much difficulty. 
\begin{assum} \label{assum:new-assumption-strong}
    The following two conditions hold:
    \begin{enumerate}[label=(\roman*)]
        \item \label{assum:new-assumption-strong-1}
        the sets of optimal solutions for \eqref{eq:hsdr} and \eqref{eq:hsdrd} are nonempty; and
        \item \label{assum:new-assumption-strong-2}
        at least one of the following two conditions holds:
        \begin{enumerate}[label=(\alph*)]
            \item \label{assum:new-assumption-strong-2-1}
                the feasible set of \eqref{eq:hsdr} is bounded; or
            \item \label{assum:new-assumption-strong-2-2}
                for \eqref{eq:hsdrd},
                the set of optimal solutions is bounded,
                and the interior of the feasible set is nonempty.
        \end{enumerate}
    \end{enumerate}
\end{assum}

We now perturb the objective function of a given QCQP
to remove the connectivity of $G$ from Theorem~\ref{thm:system-based-condition-connected}.
Let $P \in \SymMat^n$ be an $n \times n$ nonzero matrix,
and let $\varepsilon > 0$ denote the magnitude of the perturbation.
An $\varepsilon$-perturbed QCQP is described as follows:
\begin{equation}
    \label{eq:hqcqp-perturbed} \tag{$\PC^\varepsilon$}
	\begin{array}{rl}
        \min & \trans{\x} \left(Q^0 + \varepsilon P\right) \x \\
        \subto & \trans{\x} Q^p \x \le b_p, \quad p \in [m].
    \end{array}
\end{equation}
To generalize $S(\y)$ for the $\varepsilon$-perturbed QCQP,
we define
\begin{equation*}
    S(\y;\, \varepsilon)
        \coloneqq Q^0 + \varepsilon P + \sum_{p = 1}^m y_pQ^p
        = S(\y) + \varepsilon P.
\end{equation*}

\subsection{Perturbation techniques}
\label{ssec:perturbation-techniques}
Under the  condition that the feasible set of a QCQP is bounded,
Azuma et al.~\cite[Lemma 3.3]{Azuma2021} proved that the SDP relaxation is exact
if a sequence of perturbed QCQPs that satisfy the exactness condition converges to the original one.
This result was used to eliminate the requirement that
the aggregated sparsity pattern graph is connected from their main theorem.
The following lemmas are  extensions of the results in \cite{Azuma2021} under a weaker assumption.

\begin{lemma} \label{lemma:perturbation-technique-primal}
    Suppose that Assumption
    ~\ref{assum:new-assumption-strong} {\it \ref{assum:new-assumption-strong-1}} and
    \ref{assum:new-assumption-strong-2}\ref{assum:new-assumption-strong-2-1} hold.
    Let $P \neq O$ be an $n \times n$ nonzero matrix, and
    $\{\varepsilon_t\}_{t = 1}^\infty$ be a monotonically decreasing sequence
    such that $\lim_{t \to \infty} \varepsilon_t = 0$.
    If the SDP relaxation of the $\varepsilon_t$-perturbed problem
    $(\PC^{\varepsilon_t})$
    is exact for all $t = 1, 2, \ldots$,
    then the SDP relaxation of the original problem \eqref{eq:hqcqp} is also exact.
\end{lemma}
\begin{proof}
    Let $A$ and $B$ be
    the feasible sets of \eqref{eq:hqcqp} and \eqref{eq:hsdr}, respectively:
    \begin{align*}
        A \coloneqq & \left\{
            \x \in \Real^n \,\middle|\,
            \ip{Q^p}{(\x\trans{\x})} \leq b_p, \quad p = 1, \ldots, m\right\}, \\
        B \coloneqq & \left\{
            X \in \SymMat_+^n \,\middle|\,
            \ip{Q^p}{X} \leq b_p, \quad p = 1, \ldots, m\right\}.
    \end{align*}
    Note that $B$ is a compact set by the assumption.
    The intersection of $B$ and the set of rank-1 matrices
    \begin{align*}
        B_1
            &\coloneqq B \cap \left\{X \in \SymMat^n \,\middle|\, \rank(X) \leq 1\right\} \\
            &= \left\{
                X \succeq O \,\middle|\,
                \rank(X) \leq 1,\;
                \ip{Q^p}{X} \leq b_p, \; p = 1, \ldots, m\right\}
    \end{align*}
    is also a compact set since $\left\{X \in \SymMat^n \,\middle|\, \rank(X) \leq 1\right\}$ is closed.
    There exists a bijection $f: A \to B_1$ given by $f(\x) = \x\trans{\x}$,
    thus $A$ is also a compact set.
    By an argument similar to the proof of \cite[Lemma 3.3]{Azuma2021},
    we obtain the desired result.
\end{proof}

\begin{lemma} \label{lemma:perturbation-technique-dual}
    Suppose that Assumption~\ref{assum:new-assumption-strong} {\it \ref{assum:new-assumption-strong-1}} and
    \ref{assum:new-assumption-strong-2}\ref{assum:new-assumption-strong-2-2}  hold.
    Let $P \neq O$ be an $n \times n$ negative semidefinite nonzero matrix, and
    $\{\varepsilon_t\}_{t = 1}^\infty$ be a monotonically decreasing sequence
    such that $\lim_{t \to \infty} \varepsilon_t = 0$.
    If the SDP relaxation of the $\varepsilon_t$-perturbed problem
    $(\PC^{\varepsilon_t})$
    is exact for all $t = 1, 2, \ldots$,
    then the SDP relaxation of the original problem \eqref{eq:hqcqp} is also exact.
\end{lemma}
\begin{proof}
    Let $\Gamma \coloneqq \left\{\y \geq \0 \,\middle|\, S(\y) \succeq O\right\}$
    be the feasible set of \eqref{eq:hsdrd}.
    Let $(\DC_R^{\varepsilon})$ denote
    the dual of the SDP relaxation for $\varepsilon$-perturbed QCQP \eqref{eq:hqcqp-perturbed},
    and define $\Gamma(\varepsilon) \coloneqq \left\{\y \geq \0 \,\middle|\, S(\y;\, \varepsilon) \succeq O\right\}$
    as the feasible set of $(\DC_R^{\varepsilon})$.
    Since $P$ is negative semidefinite, we have $S(\y;\, \varepsilon_1) \preceq S(\y;\, \varepsilon_2)$
    for any $\y \geq \0$ and $\varepsilon_1 > \varepsilon_2 > 0$, which indicates
    a monotonic structure of the sequence $\left\{\Gamma(\varepsilon_t)\right\}_{t=1}^\infty$:
    \begin{equation*}
        \Gamma = \Gamma(0) \supseteq \cdots \supseteq \Gamma(\varepsilon_{t+1})
        \supseteq \Gamma(\varepsilon_t) \supseteq \cdots.
    \end{equation*}
    From Assumption~\ref{assum:new-assumption-strong}{\it\ref{assum:new-assumption-strong-2}\ref{assum:new-assumption-strong-2-2}},
    there exists a point $\bar{\y} \in \Gamma$ such that $S(\bar{\y}) \succ O$.
    Since each $\Gamma(\varepsilon_t)$ is a closed set and $\lim_{t \to \infty} \varepsilon_t = 0$,
    there exists an integer $T$ such that
    $S(\bar{\y}; \varepsilon_T) \succ O$. In addition, it holds that
    $S(\bar{\y}; \varepsilon_t) \succeq S(\bar{\y}; \varepsilon_T)$ for $t \ge T$.

    Let $v_t^*$ and  $B^*(\varepsilon_t)$ be the optimal value  and
    the set of the corresponding optimal solutions of $(\PC^{\varepsilon_t})$, respectively.
    From the assumptions that $(\PC)$ has a feasible point
    and $P$ is negative semidefinite,
    there is an upper bound $\bar{v}$ such that $v_t^* \le \bar{v}$ for any $t$.
    Therefore, it holds that, for any $t \ge T$,
    \begin{align*}
        B^*(\varepsilon_t)
        &= \left\{
            X \in \SymMat^n \,\middle|\,
            X \succeq O,\;
            \ip{(Q^0 + \varepsilon_t P)}{X} = v_t^*,\;
            \ip{Q^p}{X} \leq b_p \; \text{for all $p \in [m]$}
        \right\} \\
        & \subseteq \left\{
            X \in \SymMat^n \,\middle|\,
            X \succeq O,\;
            \ip{\left(Q^0 + \varepsilon_t P + \sum_{p=1}^m \bar{y}_pQ^p\right)}{X} \leq v_t^* +  \trans{\bar{\y}}\b
        \right\} \\
        &= \left\{
            X \in \SymMat^n \,\middle|\,
            X \succeq O,\;
            \ip{S(\bar{\y};\, \varepsilon_t)}{X} \leq v_t^* + \trans{\bar{\y}}\b
            \right\}, \\
        &\subseteq \left\{
            X \in \SymMat^n \,\middle|\,
            X \succeq O,\;
            \ip{S(\bar{\y};\, \varepsilon_T)}{X} \leq \bar{v} + \trans{\bar{\y}}\b
            \right\},
    \end{align*}
    which implies $\bigcup_{t=T}^\infty \; B^*(\varepsilon_t)$ is bounded
    since $S(\bar{\y};\, \varepsilon_T) \succ O$.
    With the exact SDP relaxation of the perturbed problems and strong duality,
    we can consider  $X^t \in B^*(\varepsilon_t)$, an rank-1 solution of the primal SDP relaxation,
    and $\y^t \in \Gamma(\varepsilon_t)$, an optimal solution of $(\DC_R^{\varepsilon_t})$
    satifying $X^t S(\y^t;\, \varepsilon_t) = O$.
    We define a closed set as
     \begin{equation*}
        U \coloneqq \closure\left(\bigcup_{t=T}^\infty \; B^*(\varepsilon_t)\right)
    \end{equation*}
    so that the sequence $\{X^t\}_{t=T}^\infty \subseteq U$.
    Since $\bigcup_{t=T}^\infty \; B^*(\varepsilon_t)$ is bounded, the set $U$ is a compact set.
    As the sequence has an accumulation point, we let
    $X^\mathrm{lim} \coloneqq \lim_{t \to \infty} X^t \in U$
    by taking an appropriate subsequence from $\{X^t \,|\, t \ge T\}$.
    Moreover, since $\bigcup_{t=T}^\infty \; B^*(\varepsilon_t)$ is included in the feasible set of \eqref{eq:hsdr},
    its closure $U$ is also in the same set,
    which implies that $X^\mathrm{lim}$ is an at most rank-1 feasible point of \eqref{eq:hsdr}.

    Finally, we show the optimality of $X^\mathrm{lim}$ for \eqref{eq:hsdr}.
    We assume that $\bar{X}$ is a feasible point of \eqref{eq:hsdr}
    such that $\ip{Q^0}{\bar{X}} < \ip{Q^0}{X^\mathrm{lim}}$
    and derive a contradiction.
    Since $\bigcup_{t=T}^\infty \; B^*(\varepsilon_t)$ is bounded,
    there is a sufficiently large $M$ such that
    $\| \bar{X} \| \le M$ and $\| X^{t} \| \le M$ for all $t \geq T$.
    Let $\delta = \ip{Q^0}{X^\mathrm{lim}} - \ip{Q^0}{\bar{X}} > 0$.
    Since
    $X^\mathrm{lim} = \lim_{t \to \infty} X^t $ and
    $\lim_{t \to \infty} \varepsilon_t = 0$,
    we can find $\hat{T} \ge T$ such that
    $|Q_0 \bullet (X^\mathrm{lim} - X^{\hat{T}})| \le \frac{\delta}{4}$
    and $\varepsilon_{\hat{T}} \le \frac{\delta}{8 \|P\| M }$.
    Since $\bar{X}$ and $X^{\hat{T}}$ are feasible for $(\PC^{\varepsilon_{\hat{T}}})$,
    $\frac{\bar{X} +X^{\hat{T}}}{2}$ is also feasible
    for $(\PC^{\varepsilon_{\hat{T}}})$.
    Thus, we have
    \begin{align}
    & \	\left(Q_0 + \varepsilon_{\hat{T}} P\right) \bullet \left(\frac{\bar{X} +X^{\hat{T}}}{2}\right)
        - \left(Q_0 + \varepsilon_{\hat{T}} P\right) \bullet X^{\hat{T}}  \\
    = & \ \frac{1}{2}\left(Q_0 + \varepsilon_{\hat{T}} P\right) \bullet \left(\bar{X} - X^{\hat{T}}\right) \\
    = & \ \frac{1}{2} Q_0 \bullet \left(\bar{X} - X^{\mathrm{lim}}\right)
        + \frac{1}{2} Q_0 \bullet \left(X^{\mathrm{lim}} - X^{\hat{T}}\right)
        + \frac{1}{2} \varepsilon_{\hat{T}} P \bullet \left(\bar{X} - X^{\hat{T}}\right) \\
    \le & \ \frac{1}{2} Q_0 \bullet \left(\bar{X} - X^{\mathrm{lim}}\right)
        + \frac{1}{2} \left|Q_0 \bullet \left(X^{\mathrm{lim}} - X^{\hat{T}}\right)\right|
        + \frac{1}{2} \varepsilon_{\hat{T}} \|P\| (2M) \\
    \le & -\frac{\delta}{2} +  \frac{\delta}{8} + \frac{\delta}{8}
        = - \frac{\delta}{4} < 0.
    \end{align}
     This contradicts the optimality of
    $X^{\hat{T}}$ in $(\PC^{\varepsilon_{\hat{T}}})$.
    This completes the proof.
\end{proof}

We note that the negative semidefiniteness of $P$ assumed in Lemma~\ref{lemma:perturbation-technique-dual}
is not included in Lemma~\ref{lemma:perturbation-technique-primal}.
In the subsequent discussion, we remove the assumption on the connectivity of $G$ from Theorem~\ref{thm:system-based-condition-connected}
using Lemmas~\ref{lemma:perturbation-technique-primal} and \ref{lemma:perturbation-technique-dual}.

\subsection{QCQPs with disconnected bipartite structures} \label{ssec:main-disconnected}
We present an improved version of Theorem~\ref{thm:system-based-condition-connected}
for QCQPs with disconnected aggregated sparsity pattern graphs $G$. 
\begin{theorem}
    \label{prop:system-based-condition}
    Suppose that Assumption~\ref{assum:new-assumption-strong} holds
    and that the aggregated sparsity pattern graph $G(\VC, \EC)$ is bipatite.
    Then, \eqref{eq:hsdr} is exact if, for all $(k, \ell) \in \EC$,
    the system \eqref{eq:system-nonpositive} has no solutions.
\end{theorem}
\begin{proof}
    Let $L$ denote the number of connected components of $G$,
    and choose an arbitrarily vertex $u_i$ from the connected components indexed by $i \in [L]$.
    Then, we define the edge set
    \begin{equation*}
        \FC = \bigcup_{i \in [L-1]} \left\{\left(u_i, u_{i+1}\right), \left(u_{i+1}, u_i\right)\right\}.
    \end{equation*}
    Since $\FC$ connects the $i$th and $(i+1)$th component,
    the graph $\tilde{G}(\VC, \tilde{\EC} \coloneqq \EC \cup \FC)$
    is a connected and bipartite graph.
    Let $P \in \SymMat^n$ be the negative of the Laplacian matrix of a subgraph $\hat{G}(\VC, \FC)$ of $\tilde{G}$ induced by $\FC$, i.e.,
    \begin{equation*}
        P_{ij} = \begin{cases}
            \; -\deg(i) & \quad \text{if $i = j$}, \\
            \; 1 & \quad \text{if $(i,j) \in \FC$}, \\
            \; 0 & \quad\text{otherwise},
        \end{cases}
    \end{equation*}
    where $\deg(i)$ denotes the degree of the vertex $i$ in the subgraph $\hat{G}(\VC, \FC)$.
    Since the Laplacian matrix is positive semidefinite,
    $P$ is negative semidefinite.
    By adding a perturbation $\varepsilon P$ with any $\varepsilon > 0$ into \eqref{eq:hqcqp},
    we obtain an $\varepsilon$-perturbed QCQP \eqref{eq:hqcqp-perturbed}
    whose aggregated sparsity pattern graph is $\tilde{G}(\VC, \tilde{\EC})$.

    To check the exactness of the SDP relaxation for \eqref{eq:hqcqp-perturbed}
    by Theorem~\ref{thm:system-based-condition-connected},
    it suffices to show that the following system
    \begin{equation*}
        \y \geq \0,\; S(\y;\, \varepsilon) \succeq O,\;
        S(\y;\, \varepsilon)_{k\ell} \leq 0.
    \end{equation*}
    has no solutions for all $(k, \ell) \in \tilde{\EC}$,
    where $S(\y;\, \varepsilon) \coloneqq (Q^0 + \varepsilon P) + \sum_{p \in [m]} y_p Q^p$.
    Let $\hat{\y}$ be an arbitrary vector satifying the first two constraints, i.e.,
    $\hat{\y} \geq \0$ and $S(\hat{\y};\, \varepsilon) \succeq O$.
    \begin{enumerate}[label=(\roman*)]
        \item
            If $(k, \ell) \in \FC$, then $P_{k\ell} = 1$ and
            $Q^p_{k\ell} = 0$ for any $p \in [0, m]$ by definition.
            Thus, we have
            \begin{equation*}
                S(\hat{\y};\, \varepsilon)_{k\ell}
                = \varepsilon P_{k\ell}
                > 0.
            \end{equation*}
        \item
            If $(k, \ell) \in \tilde{\EC} \setminus \FC = \EC$,
            the system \eqref{eq:system-nonpositive} with $(k, \ell)$ has no solutions,
            which implies $S(\hat{\y})_{k\ell} > 0$.
            Since $(k, \ell) \not\in \FC$, we have $P_{k\ell} = 0$.
            Hence, it follows
            \begin{equation*}
                S(\hat{\y};\, \varepsilon)_{k\ell}
                = S(\hat{\y})_{k\ell}
                > 0.
            \end{equation*}
    \end{enumerate}
    Therefore, all the systems have no solutions,
    and the SDP relaxation of \eqref{eq:hqcqp-perturbed} is exact.

    Let $\{\varepsilon_t\}_{t=1}^\infty \subseteq \Real_+$ be a monotonically decreasing sequence converging to zero,
    then the SDP relaxation of the $\varepsilon_t$-perturbed QCQP is exact as discussed above.
    By Lemmas~\ref{lemma:perturbation-technique-primal} or \ref{lemma:perturbation-technique-dual},
    the desired result follows.
\end{proof}


\subsection{Disconnected sign-definite QCQPs} \label{ssec:nonnegative-offdiagonal}
For QCQPs with the bipartite sparsity pattern and nonnegative off-diagonal elements of $Q^0, \ldots, Q^m$,
their SDP relaxation is known to be exact (see Theorem~\ref{thm:sojoudi-theorem}
\cite{Sojoudi2014exactness}).
In contrast, when we have dealt with such QCQPs in section~\ref{ssec:nonnegative-offdiagonal-connected},
the connectivity of $G$ and $Q^0_{ij} > 0$ have been assumed to derive the exactness of the SDP relaxation.
In this subsection, we eliminate these assumptions using the perturbation techniques of section~\ref{ssec:perturbation-techniques}.

\begin{coro} \label{coro:nonnegative-offdiagonal}
    Suppose that Assumption~\ref{assum:new-assumption-strong} holds, and
    suppose the aggregated sparsity pattern graph $G(\VC, \EC)$ of \eqref{eq:hqcqp}
    is bipartite.
    If $Q^p_{ij} \geq 0$ for all $(i, j) \in \EC$ and for all $p \in [0, m]$,
    then the SDP relaxation is exact.
\end{coro}
\begin{proof}
    Let $P \in \SymMat^n$ be the negative of the Laplacian matrix of $G(\VC, \EC)$, i.e.,
    \begin{equation*}
        P_{ij} = \begin{cases}
            \; -\deg(i) & \quad \text{if $i = j$}, \\
            \; 1 & \quad \text{if $(i,j) \in \EC$}, \\
            \; 0 & \quad\text{otherwise}.
        \end{cases}
    \end{equation*}
   Since the Laplacian matrix is positive semidefinite,
    $P$ is negative semidefinite.
    By adding a perturbation $\varepsilon P$ with any $\varepsilon > 0$,
    we obtain an $\varepsilon$-perturbed QCQP \eqref{eq:hqcqp-perturbed}
    whose aggregated sparsity pattern graph remains the same as the  graph $G(\VC, \EC)$.

    To determine whether the SDP relaxation is exact for this $\varepsilon$-perturbed QCQP \eqref{eq:hqcqp-perturbed},
    it suffices to check the infeasibility of the system, according to Theorem~\ref{prop:system-based-condition}:
    \begin{equation*}
        \y \geq \0,\; S(\y;\, \varepsilon) \succeq O,\;
        S(\y;\, \varepsilon)_{k\ell} \leq 0.
    \end{equation*}
    Let $\hat{\y} \geq \0$ be an arbitrary vector
    satisfying the first two constraints, i.e.,
    $\hat{\y} \geq \0$ and $S(\hat{\y};\, \varepsilon) \succeq O$.
    For every $(k, \ell) \in \EC$, since $S(\hat{\y})_{k\ell} \geq 0$ and $P_{k\ell} > 0$,
    we have
    \begin{equation*}
        S(\hat{\y};\, \varepsilon)_{k\ell}
        \geq \varepsilon P_{k\ell} > 0,
    \end{equation*}
    which implies that the system above has no solutions.
    Hence, by Theorem~\ref{prop:system-based-condition},
    the SDP relaxation of the $\varepsilon$-perturbed QCQP \eqref{eq:hqcqp-perturbed} is exact.

    Let $\{\varepsilon_t\}_{t=1}^\infty \subseteq \Real_+$ be a monotonically decreasing sequence converging to zero,
    then the SDP relaxation of the $\varepsilon$-perturbed QCQP is exact as discussed above.
    By Lemmas~\ref{lemma:perturbation-technique-primal} or \ref{lemma:perturbation-technique-dual},
    the SDP relaxation of a QCQP with nonnegative off-diagonal elements and bipartite  structures  
     is also exact.
\end{proof}

We can extend
Proposition~\ref{prop:weaker-than-sojoudi-connected}
and Corollary~\ref{coro:nonpositive-offdiagonal-connected}
using Corollary~\ref{coro:nonnegative-offdiagonal} to the following results.
\begin{prop}
    \label{prop:weaker-than-sojoudi}
    Suppose that Assumption~\ref{assum:new-assumption-strong} holds and no conditions on sparsity is considered.
    If \eqref{eq:hqcqp} satisfies the assumption of Theorem~\ref{thm:sojoudi-theorem},
    then \eqref{eq:hqcqp} also satisfies that of Corollary~\ref{coro:nonnegative-offdiagonal}.
    In addition, the exactness of its SDP relaxation
    can be proved by Theorem~\ref{prop:system-based-condition}.
\end{prop}
\begin{coro} \label{coro:nonpositive-offdiagonal}
    Under Assumption~\ref{assum:new-assumption-strong},
    the SDP relaxation of a nonpositive off-diagonal QCQP is exact.
\end{coro}

\begin{proof}
    {(Both Proposition~\ref{prop:weaker-than-sojoudi} and Corollary~\ref{coro:nonpositive-offdiagonal})}
    It is easy to check that
    the aggregated sparsity pattern graph of
    \eqref{eq:decomposed-hqcqp} generated by the given problem is bipartite
    by the arguments similar to the proof of Proposition~\ref{prop:weaker-than-sojoudi-connected}.
    Therefore, \eqref{eq:decomposed-hqcqp} satisfies the assumption of Corollary~\ref{coro:nonnegative-offdiagonal}.
\end{proof}

\section{Numerical experiments} \label{sec:example}

We investigate analytical and computational aspects of the conditions in  
 Theorem~\ref{thm:system-based-condition-connected}
with two QCQP instances below. 
The first QCQP consists of $2 \times 2$ data matrices. 
We  show the exactness of its SDP relaxation
by checking the feasibility systems in Theorem~\ref{thm:system-based-condition-connected} without   SDP solvers.
Next, Example~\ref{examp:cycle-graph-4-vertices} is considered for the second QCQP.
As the size $n$ of the second QCQP is 4, it is difficult to handle
the positive semidefinite constraint $S(\y) \succeq O$ without numerical computation. 
We present a numerical method for testing the exactness of the SDP relaxation with a computational solver.

We also detail the difference between our results and the existing results
using these two QCQP instances.
As discussed in section~\ref{ssec:comparison},
if the aggregated sparsity pattern graph is bipartite,
then Theorem~\ref{thm:system-based-condition-connected} covers a wider class of QCQPs than
those by Theorem~\ref{thm:sojoudi-theorem} in~\cite{Sojoudi2014exactness}
under the connectivity and the elementwise condition on $Q^0$.
Theorem~\ref{thm:system-based-condition-connected} has been generalized in section~\ref{sec:perturbation} to Theorem~\ref{prop:system-based-condition},
and this theorem covers a wider class of QOCPs without the connectivity condition.

For numerical experiments,
JuMP~\cite{Dunning2017} was used  with the solver MOSEK~\cite{mosek}
and  SDPs were solved with tolerance $1.0 \times 10^{-8}$. 
All numerical results are shown with four significant digits.

\subsection{A QCQP instance with $n=2$} \label{ssec:analytical-example}
\begin{example}
    \label{examp:small-example}
    Consider the QCQP \eqref{eq:hqcqp} with
    \begin{align*}
        & n = 2, \quad m = 1, \quad \b = \begin{bmatrix} 1 \end{bmatrix}, \\
        & Q^0 = \begin{bmatrix} -3 & -1 \\ -1 & -2 \end{bmatrix}, \quad
          Q^1 = \begin{bmatrix}  3 &  4 \\  4 &  6 \end{bmatrix}.
    \end{align*}
\end{example}

We first verify whether the problem satisfies the assumption of Theorem~\ref{thm:system-based-condition-connected}.
The aggregated sparsity pattern graph $G$ is bipartite and connected
as it has only two vertices and $Q^0_{12} \neq 0$.
Since $Q^1$ is positive definite,
the problem satisfies Assumption~\ref{assum:previous-assumption}{\it \ref{assum:previous-assumption-1}}.
By the discussion in Remark~\ref{rema:comparison-assumption},
it also satisfies Assumption~\ref{assum:new-assumption}.
It only remains to show that the system
\begin{equation*}
    y_1 \geq 0, \quad
    \hat{S}(y_1) \coloneqq \begin{bmatrix} -3 & -1 \\ -1 & -2 \end{bmatrix} +
        y_1 \begin{bmatrix}  3 &  4 \\  4 &  6 \end{bmatrix} \succeq O, \quad
    -1 + 4y_1 \leq 0
\end{equation*}
has no solutions.
By definition,
$\hat{S}(y_1) \succeq O$ holds if and only if all the principal minors of $\hat{S}(y_1)$ are nonnegative,
or equivalently, $-3 + 3y_1 \geq 0$, $-2 + 6y_1\ge 0$, and $2y_1^2 - 16y_1 + 5 \geq 0$.
Hence, if $y_1 \geq 4 + 3\sqrt{6}/2 \simeq 7.674$, then
the first two inequalities of the system are satisfied.
Since $-1 + 4y_1 \geq -1 + 4(4 + 3\sqrt{6}/2) = 15 + 6\sqrt{6} > 0$,
the last inequality does not hold for such  $y_1$.
The problem therefore admits the exact SDP relaxation.

Actually, we numerically obtained an optimal solution of the above QCQP in
Example~\ref{examp:small-example} and its SDP relaxation
		as $\x^* \simeq [1.731; -1.167]$ and
		$X^* \simeq [2.997, -2.021; -2.021, 1.362]$, respectively.
	From $\trans{(\x^*)}Q^0 \x^* - Q^0 \bullet X^* \simeq 5.379 \times 10^{-10}$,
		we see numerically that the SDP relaxation provided the exact optimal value. 

Since $G$ is clearly a forest (no cycles),
we can also apply Proposition~\ref{prop:forest-results} in~\cite{Azuma2021}.
From the  discussion above,
the system \eqref{eq:system-zero} has no solutions for $(k, \ell) = (1, 2)$
and Assumption~\ref{assum:previous-assumption}{\it \ref{assum:previous-assumption-1}} is satisfied.
By taking $\hat{X} = [0.1 \  \  0; 0 \ \ 0.1] \succ O$,
we know $\ip{Q^1}{\hat{X}} = 0.9 \leq 1 = b_1$.
Hence, the exactness of the SDP relaxation can be proved by Proposition~\ref{prop:forest-results}.
We mention that this result cannot be obtained by
Theorem~\ref{thm:sojoudi-theorem} in~\cite{Sojoudi2014exactness}. 
Since $Q^0_{12} = -1$ and $Q^1_{12} = 4$,
the edge sign $\sigma_{12}$ of the edge $(1, 2)$ must be zero by definition,
contradicting \eqref{eq:sign-constraint-sign-definite}.

\subsection{Example~\ref{examp:cycle-graph-4-vertices}} \label{ssec:computational-example}

We computed an optimal solution of Example~\ref{examp:cycle-graph-4-vertices} and that of its SDP relaxation as
\begin{equation*}
    x^* \simeq \begin{bmatrix}
        7.818 \\ -8.331 \\ 1.721 \\ -7.019
    \end{bmatrix}\ \text{and} \
    X^* \simeq \begin{bmatrix}
         61.12 & -65.13 &  13.45 & -54.87 \\
        -65.13 &  69.41 & -14.34 &  58.48 \\
         13.45 & -14.34 &  2.961 & -12.08 \\
        -54.87 &  58.48 & -12.08 &  49.27
    \end{bmatrix} \in \SymMat^4,
\end{equation*}
respectively.
From $\trans{(\x^*)}Q^0 \x^* - Q^0 \bullet X^* \simeq 7.676 \times 10^{-8}$,
we see numerically that the SDP relaxation resulted in the exact optimal value.

The aggregated sparsity pattern graph $G(\VC, \EC)$ is
a cycle graph with 4 vertices (\autoref{fig:example-aggregated-sparsity}).
We first see whether
it satisfies the assumption of Theorem~\ref{thm:system-based-condition-connected}.
We compute  $3Q_1 + 4Q_2$ as
\begin{equation*}
    3 \begin{bmatrix}
     5 &  2 & 0 &  1 \\  2 & -1 &  3 &  0 \\
     0 &  3 & 3 & -1 \\  1 &  0 & -1 &  4 \end{bmatrix} +
    4 \begin{bmatrix}
    -1 &  1 & 0 &  0 \\  1 &  4 & -1 &  0 \\
     0 & -1 & 6 &  1 \\  0 &  0 &  1 & -2 \end{bmatrix} =
    \begin{bmatrix}
    11 & 10 &  0 &  3 \\ 10 & 13 &  5 &  0 \\
     0 &  5 & 33 &  1 \\  3 &  0 &  1 &  4 \end{bmatrix},
\end{equation*}
and its minimum eigenvalue is approximately $0.1577$.
Thus, there exists $\bar{\y} \geq 0$ such that $\bar{y}_1 Q_1 + \bar{y}_2 Q_2 \succ O$,
e.g., $\bar{\y} = [3; 4]$.
As mentioned in Remark~\ref{rema:comparison-assumption},
it follows that the second problem satisfies Assumption~\ref{assum:new-assumption}.
To show the exactness of the SDP relaxation for the problem, 
it only remains to show that
the systems \eqref{eq:system-nonpositive} for all $(k, \ell) \in \EC$ has no solutions.
Using an SDP solver on a computer, we could observe that there is no solution for the system.
%
Indeed, for every $(k, \ell) \in \EC$, the SDP
\begin{equation} \label{eq:optimal-value-systems-example}
    \begin{array}{rl}
        \mu^* =
        \min & S(\y)_{k\ell} \\
        \subto & \y \geq \0, \; S(\y) \succeq O,
    \end{array}
\end{equation}
returns the optimal values shown in \autoref{tab:optimal-value-systems-example},
which implies that no solution exists for \eqref{eq:system-nonpositive}
since $S(\y)_{k\ell}$ cannot attain a nonpositive value.
Therefore,
the SDP relaxation of Example~\ref{examp:cycle-graph-4-vertices} is exact by
Theorem~\ref{thm:system-based-condition-connected}.
\begin{table}[b]
	\caption{Optimal values of \eqref{eq:optimal-value-systems-example} for each $(k, \ell)$}
	\label{tab:optimal-value-systems-example}
	\centering
	\begin{tabular}{c|cccc}
		$(k, \ell)$ & $(1, 2)$ & $(2, 3)$ & $(1, 4)$ & $(3, 4)$ \\ \hline
		$\mu^*$ & 18.58 & 12.84 & 8.897 & 0.3215
	\end{tabular}
\end{table}

With Theorem~\ref{thm:sojoudi-theorem} in~\cite{Sojoudi2014exactness},
it is not possible to show the exactness of the SDP relaxation. 
The edge sign $\sigma_{12}$ for $(1, 2)$th element is $0$ by definition.
Since the cycle basis of $\GC$ is only $\CC_1 = \GC$,
the left-hand side of \eqref{eq:sign-constraint-simple-cycle} is
$\sigma_{12}\sigma_{23}\sigma_{34}\sigma_{41} = 0$.
However, its right-hand side only takes $-1$ or $+1$.
This implies that Theorem~\ref{thm:sojoudi-theorem} cannot be applied to Example~\ref{examp:cycle-graph-4-vertices}.

\section{Concluding remarks} \label{sec:conclution}
We have proposed sufficient conditions for
the exact SDP relaxation of QCQPs whose
aggregated sparsity pattern graph can be represented by bipartite graphs.
Since these conditions consist of at most $n^2/4$ SDP systems,
the exactness can be investigated in polynomial time.
The derivation of the conditions is based on 
the rank of optimal solutions $\y$ of the dual SDP relaxation under strong duality.
More precisely, a QCQP admits the exact SDP relaxation
if the lower bound of the rank of $S(\y)$ is $n-1$.
For the lower bound,
we have used the fact that
any nonnegative matrix $M \succeq O$ with bipartite sparsity pattern is of at least rank $n - 1$
if it satisfies $M\1 > \0$.

Using results from the recent paper~\cite{kim2021strong},
the sufficient conditions have been considered under weaker assumptions than
those in \cite{Azuma2021}.
That is,  the sparsity of bipartite graphs includes that of tree and forest graphs,
therefore, the proposed conditions
can serve for a wider class of QCQPs 
than those in \cite{Azuma2021}.
We have also shown in Proposition~\ref{prop:weaker-than-sojoudi} that
one can determine the exactness for all the problems
which satisfy the condition considered in Theorem~\ref{thm:sojoudi-theorem} (\cite{Sojoudi2014exactness}).

For our future work,
 sufficient conditions for the exactness  of
a wider class of QCQPs than those with bipartite structures will be investigated.
Furthermore, examining our conditions
to analyze the exact SDP relaxation of QCQPs
transformed from polynomial optimization would be an interesting subject.

\vspace{0.5cm}

\noindent
{\bf Acknowledgements.}
The authors would like to thank Prof. Ram Vasudevan and Mr. Jinsun Liu  for pointing out that there exists 
no edge $(i,i+n)$ in the objective function  in the proof of Proposition 3.8
of the original version.

\bibliographystyle{abbrv} 

\bibliography{./reference}

\begin{thebibliography}{10}

\bibitem{Anton2014}
H.~Anton and C.~Rorres.
\newblock {\em Elementary Linear Algebra: Applications Version}.
\newblock John Wiley \& Sons Inc., USA, 11th ed. edition, 2014.

\bibitem{Argue2020necessary}
C.~J. Argue, F.~K{\dotlessi}l{\dotlessi}n\chige{c}-Karzan, and A.~L. Wang.
\newblock Necessary and sufficient conditions for rank-one generated cones.
\newblock arXiv:2007.07433, 2020.

\bibitem{Azuma2021}
G.~Azuma, M.~Fukuda, S.~Kim, and M.~Yamashita.
\newblock Exact {{SDP}} relaxations of quadratically constrained quadratic
  programs with forest structures.
\newblock {\em Journal of Global Optimization}, 82(2):243--262, 2022.

\bibitem{BISWAS2004}
P.~Biswas and Y.~Ye.
\newblock Semidefinite programming for ad hoc wireless sensor network
  localization.
\newblock In {\em Proceedings of the Third International Symposium on
  Information Processing in Sensor Networks}, pages 46--54, New York, 2004.
  ACM.

\bibitem{Burer2019}
S.~Burer and Y.~Ye.
\newblock Exact semidefinite formulations for a class of (random and
  non-random) nonconvex quadratic programs.
\newblock {\em Mathematical Programming}, 181(1):1--17, 2020.

\bibitem{Dunning2017}
I.~Dunning, J.~Huchette, and M.~Lubin.
\newblock Jump: A modeling language for mathematical optimization.
\newblock {\em SIAM Review}, 59(2):295--320, 2017.

\bibitem{Geomans1995}
M.~X. Goemans and D.~P. Williamson.
\newblock Improved approximation algorithms for maximum cut and satisfiability
  problems using semidefinite programming.
\newblock {\em Journal of the ACM}, 42(6):1115--1145, 1995.

\bibitem{grone1992nonchordal}
R.~Grone, R.~Loewy, and S.~Pierce.
\newblock Nonchordal positive semidefinite stochastic matrices.
\newblock {\em Linear and Multilinear Algebra}, 32(2):107--113, 1992.

\bibitem{Hsia2013}
Y.~Hsia and R.-L. Sheu.
\newblock Trust region subproblem with a fixed number of additional linear
  inequality constraints has polynomial complexity.
\newblock arXiv:1312.1398, 2013.

\bibitem{Jeyakumar2014}
V.~Jeyakumar and G.~Y. Li.
\newblock Trust-region problems with linear inequality constraints: Exact sdp
  relaxation, global optimality and robust optimization.
\newblock {\em Mathematical Programming}, 147(1-2):171--206, 2014.

\bibitem{kim2003exact}
S.~Kim and M.~Kojima.
\newblock Exact solutions of some nonconvex quadratic optimization problems via
  {SDP} and {SOCP} relaxations.
\newblock {\em Computational Optimization and Applications}, 26(2):143--154,
  2003.

\bibitem{kim2021strong}
S.~Kim and M.~Kojima.
\newblock Strong duality of a conic optimization problem with a single
  hyperplane and two cone constraints.
\newblock arXiv:2111.03251v2, 2021.

\bibitem{KIM2009}
S.~Kim, M.~Kojima, and H.~Waki.
\newblock Exploiting sparsity in {SDP} relaxation for sensor network
  localization.
\newblock {\em SIAM {J}ournal on {O}ptimization}, 20(1):192--215, 2009.

\bibitem{kimizuka2019solving}
M.~Kimizuka, S.~Kim, and M.~Yamashita.
\newblock Solving pooling problems with time discretization by {LP} and {SOCP}
  relaxations and rescheduling methods.
\newblock {\em Journal of Global Optimization}, 75(3):631--654, 2019.

\bibitem{Lavaei2012}
J.~Lavaei and S.~H. Low.
\newblock Zero duality gap in optimal power flow problem.
\newblock {\em IEEE Transactions on Power Systems}, 27(1):92--107, 2012.

\bibitem{Locatelli2016}
M.~Locatelli.
\newblock Exactness conditions for an {SDP} relaxation of the extended trust
  region problem.
\newblock {\em Optimization Letters}, 10(6):1141--1151, 2016.

\bibitem{mosek}
{MOSEK ApS}.
\newblock Mosek/mosek.jl: Interface to the {M}osek solver in {J}ulia, 2022.
\newblock \url{https://github.com/MOSEK/Mosek.jl} (accessed on April 9, 2022).

\bibitem{Polik2007}
I.~P{\'o}lik and T.~Terlaky.
\newblock A survey of the {S}-lemma.
\newblock {\em SIAM Review}, 49(3):371--418, 2007.

\bibitem{PRendl09}
J.~Povh and F.~Rendl.
\newblock Copositive and semidefinite relaxations of the quadratic assignment
  problem.
\newblock {\em Discrete Optimization}, 6(3):231--241, 2009.

\bibitem{sheen2020exploiting}
H.~Sheen and M.~Yamashita.
\newblock Exploiting aggregate sparsity in second-order cone relaxations for
  quadratic constrained quadratic programming problems.
\newblock {\em Optimization Methods and Software}, pages 1--19, 2020.

\bibitem{SO2007}
A.~M. So and Y.~Ye.
\newblock Theory of semidefinite programming for sensor network localization.
\newblock {\em Mathematical Programming}, 109(2--3):367--384, 2007.

\bibitem{Sojoudi2014exactness}
S.~Sojoudi and J.~Lavaei.
\newblock Exactness of semidefinite relaxations for nonlinear optimization
  problems with underlying graph structure.
\newblock {\em SIAM Journal on Optimization}, 24(4):1746--1778, 2014.

\bibitem{Wang2021geometric}
A.~L. Wang and F.~K{\dotlessi}l{\dotlessi}n\chige{c}-Karzan.
\newblock A geometric view of {SDP} exactness in {QCQPs} and its applications.
\newblock arXiv:2011.07155v3, 2021.

\bibitem{Wang2021tightness}
A.~L. Wang and F.~K{\dotlessi}l{\dotlessi}n\chige{c}-Karzan.
\newblock On the tightness of {SDP} relaxations of {QCQPs}.
\newblock {\em Mathematical Programming}, 2021.

\bibitem{Wang2015}
S.~Wang and Y.~Xia.
\newblock Strong duality for generalized trust region subproblem: S-lemma with
  interval bounds.
\newblock {\em Optimization Letters}, 9(6):1063--1073, 2015.

\bibitem{Yakubovich1971}
V.~A. Yakubovich.
\newblock S-procedure in nonlinear control theory.
\newblock {\em Vestnik Leningrad University Mathematics}, 1:62--77, 1971.

\bibitem{ZHAO1998}
Q.~Zhao, S.~Karisch, F.~Rendl, and H.~Wolkowicz.
\newblock Semidefinite programming relaxations for the quadratic assignment
  problem.
\newblock {\em Journal of Combinatorial Optimization}, 2(1):71--109, 1998.

\bibitem{Zhou2019}
F.~Zhou, Y.~Chen, and S.~H. Low.
\newblock Sufficient conditions for exact semidefinite relaxation of optimal
  power flow in unbalanced multiphase radial networks.
\newblock {\em IEEE 58th {C}onference on {D}ecision and {C}ontrol (CDC)},
  58:6227--6233, 2019.

\end{thebibliography}

\end{document}